  \documentclass[10pt]{amsart}
\usepackage{amsmath}
\usepackage{amsxtra}
\usepackage{amscd}
\usepackage{amsthm}
\usepackage{amsfonts}
\usepackage{amssymb}
\usepackage{eucal}

\theoremstyle{definition}

\theoremstyle{remark}

\newcommand{\thmref}[1]{Theorem~\ref{#1}}
\newcommand{\secref}[1]{\S\ref{#1}}
\newcommand{\lemref}[1]{Lemma~\ref{#1}}
\newcommand{\propref}[1]{Proposition~\ref{#1}}

\newcommand{\remref}[1]{Remark~\ref{#1}}
\newcommand{\exref}[1]{Example~\ref{#1}}

\newcommand{\nc}{\newcommand}
\nc{\renc}{\renewcommand}
 \def\tsec#1{\medskip {\noindent \bf #1}.}
 \nc{\ssec}{\subsection}
\nc{\sssec}{\subsubsection}
\nc{\on}{\operatorname}
\nc{\remm}[1]{\<{remark} \ \lbl{#1} \>{remark}}

\nc\ol{\overline}
\nc\wt{\widetilde}
\nc\wh{\widehat}
\nc\tboxtimes{\wt{\boxtimes}}

\nc{\Aa}{{\mathbb{A}}}
 \nc{\Gg}{{\mathbb{G}}}  
\nc{\ba}{{  \bar{a}}}
\nc{\bb}{  \bar{b}}
 \nc{\Nn}{{\mathbb{N}}}
\nc{\Pp}{{\mathbb{P}}}
\nc{\Rr}{{\mathbb{R}}}
\nc{\BV}{{\mathbb{V}}}
\nc{\BW}{{\mathbb{W}}}
\nc{\Zz}{{\mathbb{Z}}}
\nc{\Qq}{{\mathbb{Q}}}
\nc{\Ss}{{\mathbb{S}}}
\nc{\Cc}{{\mathbb{C}}}
\nc{\Ff}{{\mathbb{F}}}
\nc{\Oo}{{\mathcal{O}}}
\nc{\Mm}{{\mathcal{M}}}
\nc{\dU}{{\overset{\bullet}{\bigcup}}{}}
\def\dl{\delta}
\def\du{\dU}
\nc{\dual}[1]{{\overset{\vee}{#1}{}}}

\nc{\cM}{{\check{\mathcal M}}{}}
 \nc{\oM}{{\overset{\circ}{\mathcal M}}{}}

\nc{\fD}{{\mathfrak{D}}}
\nc{\fa}{{\mathfrak{a}}}
\nc{\fb}{{\mathfrak{b}}}
\nc{\fg}{{\mathfrak{g}}}
\nc{\fgl}{{\mathfrak{gl}}}
\nc{\fh}{{\mathfrak{h}}}
\nc{\fj}{{\mathfrak{j}}}
\nc{\fm}{{\mathfrak{m}}}
\nc{\fn}{{\mathfrak{n}}}
\nc{\fu}{{\mathfrak{u}}}
\nc{\fp}{{\mathfrak{p}}}
\nc{\fr}{{\mathfrak{r}}}
\nc{\fs}{{\mathfrak{s}}}
\nc{\fsl}{{\mathfrak{sl}}}
\nc{\hsl}{{\widehat{\mathfrak{sl}}}}
\nc{\hgl}{{\widehat{\mathfrak{gl}}}}
\nc{\hg}{{\widehat{\mathfrak{g}}}}
\nc{\chg}{{\widehat{\mathfrak{g}}}{}^\vee}
\nc{\hn}{{\widehat{\mathfrak{n}}}}
\nc{\chn}{{\widehat{\mathfrak{n}}}{}^\vee}

\nc{\tT}{{\mathfrak{T}}} 
\def\tTt{\sA}
\nc{\fF}{{\mathfrak{F}}}

\nc{\bc}{{\mathbf{c}}}
\nc{\bd}{{\mathbf{d}}}
\nc{\be}{{\mathbf{e}}}
\nc{\bj}{{\mathbf{j}}}
\nc{\bn}{{\mathbf{n}}}
\nc{\bp}{{\mathbf{p}}}
\nc{\bq}{{\mathbf{q}}}
\nc{\bF}{{\mathbf{F}}}
\nc{\bu}{{\mathbf{u}}}
\nc{\bv}{{\mathbf{v}}}
\nc{\bx}{{\mathbf{x}}}
\nc{\bs}{{\mathbf{s}}}
\nc{\by}{{\mathbf{y}}}
\nc{\bw}{{\mathbf{w}}}
\nc{\bA}{{\mathbf{A}}}
\nc{\bK}{{\mathbf{K}}}
\nc{\bI}{{\mathbf{I}}}
\nc{\bB}{{\mathbf{B}}}
\nc{\bG}{{\mathbf{G}}}
\nc{\bC}{{\mathbf{C}}}
\nc{\bD}{{\mathbf{D}}}
\nc{\bP}{{\mathbf{P}}}
\nc{\bH}{{\mathbf{H}}}
\nc{\bM}{{\mathbf{M}}}
\nc{\bN}{{\mathbf{N}}}
\nc{\bV}{{\mathbf{V}}}
\nc{\bU}{{\mathbf{U}}}
\nc{\bL}{{\mathbf{L}}}
\nc{\bT}{{\mathbf{T}}}
\nc{\bW}{{\mathbf{W}}}
\nc{\bX}{{\mathbf{X}}}
\nc{\bY}{{\mathbf{Y}}}
\nc{\bZ}{{\mathbf{Z}}}
\nc{\bS}{{\mathbf{S}}}

\nc{\sA}{{\mathsf{A}}}
\nc{\sB}{{\mathsf{B}}}
\nc{\sC}{{\mathsf{C}}}
\nc{\sD}{{\mathsf{D}}}
\nc{\sF}{{\mathsf{F}}}
\nc{\sG}{{\mathsf{G}}}
\nc{\sK}{{\mathsf{K}}}

\nc{\sO}{{\mathsf{O}}}
\nc{\sQ}{{\mathsf{Q}}}
\nc{\sP}{{\mathsf{P}}}
\nc{\sZ}{{\mathsf{Z}}}
\nc{\sfp}{{\mathsf{p}}}
\nc{\sr}{{\mathsf{r}}}
\nc{\sg}{{\mathsf{g}}}
\nc{\sff}{{\mathsf{f}}}
\nc{\sfb}{{\mathsf{b}}}
\nc{\sfc}{{\mathsf{c}}}
\nc{\sd}{{\mathsf{d}}}

\nc{\tA}{{\widetilde{{A}}}}
\nc{\tD}{{\widetilde{{D}}}}
\nc{\tH}{{\widetilde{{H}}}}
\nc{\tB}{{\widetilde{{B}}}}
\nc{\tM}{{\widetilde{{M}}}}
\nc{\tg}{{\widetilde{\mathfrak{g}}}}
\nc{\tG}{{\widetilde{G}}}
\nc{\TM}{{\widetilde{\mathbb{M}}}{}}
\nc{\tO}{{\widetilde{\mathsf{O}}}{}}
\nc{\tU}{{\widetilde{\mathfrak{U}}}{}}
\nc{\TZ}{{\tilde{Z}}}
\nc{\tx}{{\tilde{x}}}
\nc{\tbv}{{\tilde{\bv}}}
\nc{\tfP}{{\widetilde{\mathfrak{P}}}{}}
\nc{\tz}{{\tilde{\zeta}}}
\nc{\tmu}{{\tilde{\mu}}}

  \nc{\Sym}{{\mathop{\operatorname{\rm Sym}}}}
   \nc{\Aut}{{\mathop{\operatorname{\rm Aut}}}}
 \nc{\Spec}{{\mathop{\operatorname{\rm Spec}}}}
\nc{\Ker}{{\mathop{\operatorname{\rm Ker}}}}
 \nc{\dom}{{\mathop{\operatorname{\rm dom}}}}
\nc{\End}{{\mathop{\operatorname{\rm End}}}}
 \nc{\Hom}{\on{\Hom}}
 \nc{\GL}{{\mathop{\operatorname{\rm GL}}}}
 \nc{\Id}{{\mathop{\operatorname{\rm Id}}}}
 \nc{\rk}{{\mathop{\operatorname{\rm rk}}}}
\nc{\irk}{{\mathop{\operatorname{\rm i-rk}}}}
 \nc{\length}{{\mathop{\operatorname{\rm length}}}}
\nc{\supp}{{\mathop{\operatorname{\rm supp}}}}
\nc{\val}{{\rm val}}
\nc{\res}{{\mathop{\operatorname{\rm res}}}}

\def\tensor{{\otimes}}
\def\meet{\cap}
\def\union{\cup}

\def\si{\sigma}

\def\<{\begin}
 \def\>{\end}

\nc{\tV}{{\widetilde{{V}}}}
\nc{\hb}[1]{\hbox{#1}}

\nc{\sM}{{\setminus}}
\renc{\-}{\setminus}
 
\nc{\seq}[1]{\stackrel{#1}{\sim}}
   
\def\inv{{^{-1}}}
\def\claim{{\bf Claim \ }}
\def\beq#1{   \begin{equation} \label{#1}   }
\def\eeq{\end{equation}   }

\def\Uu{\mathbb U} 

 \def\Vv{\mathbb V}

  \def\Vvb{ { \Vv ^b}}
  \def\St\Vvb
 
\def\iso{\simeq}

\def\prf{\begin{proof}}
\def\eprf{\end{proof} }

\def\acl{\mathop{\rm acl}\nolimits}

\def\lbl#1{  \label{#1}  }

 \author{Ehud Hrushovski}
 
\address{\newline Institute of Mathematics, the Hebrew
 University of Jerusalem, Givat Ram, Jerusalem, 91904, Israel.} 
 
 \email{ehud@math.huji.ac.il}

 \thanks{Thanks to Yad Hanadiv, and to the ISF,   grant 1048/07.  Thanks to the referee for many comments.  MSC 03c99, 18B40. }

\def\th{{\tilde{h}}}
\usepackage{color}

\title{Groupoids, imaginaries and internal covers}

\begin{document}

 \maketitle
 \<{abstract}{ Let $T$ be a first-order theory.  A  correspondence is established  between internal covers of 
models of $T$ and definable groupoids within $T$.  We also consider amalgamations
of independent diagrams of algebraically closed substructures, and find strong relation between:   covers,  uniqueness for 3-amalgamation, existence of 4-amalgamation,  imaginaries of $T^\si$, and definable groupoids.   As a corollary, we describe
the imaginary elements of families of finite-dimensional vector spaces over pseudo-finite fields.} \end{abstract}
 \def\tf{{\tilde{f}}}
 \def\tri{{\mathsf B}}
 \def\tres{{\mathbf{t}}}
 \def\tTs{\widetilde{T_\si}}
 \def\fG{{\mathcal G}}
 \def\fC{{\mathcal C}}
 \def\tC{{   \widetilde{\fC}   }}
\def\fGV{{\fG}_{\Vv}}
  \nc{\Ob}{{\mathop{\operatorname{\rm Ob}}}}
    \nc{\Mor}{{\mathop{\operatorname{\rm Mor}}}}
        \nc{\dcl}{{\mathop{\operatorname{\rm dcl}}}}
\def\ObG{\Ob {\fG}}  
\def\ObGV{{\Ob}{\fG_{\Vv}}}
\def\MorG{{\Mor}{\fG}}
\renc\k{\mathbf{k}}

\def\G{\Gamma}
\def\RV{RV}  \def \rv{rv}

\def\iso{\cong}

The questions this manuscript addresses arose in the course of an investigation of the
imaginary sorts in  ultraproducts of $p$-adic fields.  
These were shown to  be understandable given 
 the imaginary sorts of  certain finite-dimensional vector spaces over the residue field.  The residue field is pseudo-finite, and the imaginary elements there were previously studied,
 and shown in fact to be eliminable over an appropriate base. 
   It remains therefore to  
describe the imaginaries of  finite-dimensional vector spaces over a field $F$,
given those of $F$.
I expected this step to be rather easy; but it turned out to become easy only
after  a number of issues, of interest in themselves, are made clear.  

Let $T$ be a first-order theory.  A  correspondence is established  between internal covers of 
models of $T$ and definable groupoids within $T$.  
Internal covers were recognized as central in the study of totally categorical structures, but nevertheless remained mysterious; 
it was not clear how to describe the possible $T'$ from the point of view of $T$.  We give an account of this here, in terms of groupoids in place
of equivalence relations.   This description permits the view of the cover as a generalized imaginary sort.

This seems to be a useful language even for finite covers, though 
there the situation is rather  well-understood, cf.  \cite{ev}.   
   We concentrate on finite generalized
imaginaries, and describe a  a connection  between elimination of imaginaries and  higher amalgamation principles within the algebraic closure of an independent $n$-tuple.  
   The familiar imaginaries of $T^{eq}$ correspond to 3-amalgamation, as was understood for some time for stable and simple theories, 
and finite generalized imaginaries correspond to   4-amalgamation. 
This brings out ideas present in some form in \cite{ch-h}, \cite{bouh}, \cite{evh}, \cite{ev}.  
In particular, 4-amalgamation always holds for stable theory $T$, if ``algebraic closure''
is taken to include generalized imaginaries.  We also relate uniqueness of $n$-amalgamation to
existence of $n+1$-amalgamation; using ``all'' finite imaginaries (not necessarily arising from groupoids) we show that   $n$-amalgamation exists and is unique for all $n$.    

Adding an automorphism to the language to obtain a Robinson theory $T^\si$ has the effect of shifting the amalgamation 
dimension by one;  $n$-amalgamation in the expanded language corresponds
to   $n+1$-amalgamation for $T$.  Thus ordinary imaginaries of $T^\si$ can be understood, given generalized imaginaries of $T$.  

We thus find a strong relation between four things:   covers, failure of uniqueness for 3-amalgamation, imaginaries of $T^\si$,
and definable groupoids.   A clear continuation to $n=4$ would be interesting.  
 
Returning to the original motivation, we use these ideas to determine the imaginaries for systems of finite-dimensional
vector spaces over fields, and especially over pseudo-finite fields (\thmref{pfeli}).

  Thanks to Levon Haykazyan and Rahim Moosa for correcting the hypothesis of \thmref{im-gr}, and
 to Christian d'Elbée for his useful comments on the proof of \remref{stableEIr}.
 
\<{section}{Preliminaries}  \lbl{prelim}

Let $T$ be a first-order theory, with universal domain $\Uu$. 
   \footnote{More generally we can work
 with a ``Robinson theory'', 
a universal theory with the amalgamation property for substructure; one then works
with substructures of a universal domain, and takes ``definable'' to mean:  quantifier-free definable.  This was one of the ``contexts" of   \cite{lazy}; I dubbed it ``Robinson" when unaware of this reference, and the name stuck.}
$Def(\Uu)$ is the category of $\Uu$-definable sets (with parameters)
and maps between them.   

 Let $A,B$ be  small subsets of $\Uu$.  For each $b \in B$, we provide a new constant symbol
 $c_b$; and for each $a \in A$, a new variable $x_a$.  We write $tp(A/B)$ for the 
 set of all formulas with these new variables and constants, true in $\Uu$ under the eponymous
 interpretation of constant symbols and assignment of variables.  This is useful in
 expressions such as $tp(A/B) \models tp(A/B')$.  
 
An $\infty$-definable set is the solution set of a partial type (of bounded size; say bounded by the cardinality of the language.)  
 Morphisms between $\infty$ -definable sets are still induced by ordinary definable maps.  If the partial
type is allowed to have infinitely many variables, the set is called
$\star$-definable instead.  $\star$-definable sets can also be viewed
as projective systems of definable sets and maps.  

Dually, a $\bigvee$-definable set is the complement of an $\infty$-definable set.

When we say a set $P$ is  definable, we mean:   without parameters. 
If we wish to speak about a set definable with parameters $a$, we will
exhibit these parameters in the notation: $P_a$.

We will often consider two languages $L \subset L'$.  The language $L'$ may have more sorts than $L$.  Let $T'$ be a complete theory for $L'$, $T=T'|L$.
We say $T$ is {\em embedded} if any relation of $L$ is $T'$-equivalent to a formula of $L'$.   We say the sorts of $L$ are {\em stably embedded} if
in any model $M' \models T'$, any $M'$-definable subset of $S_1 \times \times S_k$ (where the $S_i$ are $L$-sorts) is also definable with parameters from $\union S_i(M')$.
This basic notion has various equivalent forms, see appendix to \cite{CH} and also \cite{ah-z}.

Let $D$ be a definable set of $L'$.  We say $D$ is {\em internal} to $L$ if in some (or any) model $M'$ of $T'$, there exist sorts $S_1,\ldots,S_k$ of $L$
and an $M$-definable   map $f$ whose domain is a subset of $S_1 \times \times S_k$, and whose image is $D$.   See \cite{bedlewo}, appendix, where it is shown that internality is associated
with definable automorphism groups; indeed, assuming $T$ is embedded and stably embedded in $T'$,  and $L' \setminus L$ is finite for simplicity, and letting $M$ denote the $L'$-sorts of $M'$, there exists a definable
group $G$ such that $G(M')$ can be identified with  $Aut(D(M')/M)$.  $G$ is called the {\em liaison group}, a term due to Poizat.    It is also shown in \cite{bedlewo} that $G$ is $M$-isomorphic to an $M$-definable group.  In \secref{liaison} we will prove a    more precise,  parameter-free version, using the notion of a {\em definable groupoid}.

We can  immediately introduce   one of the main notions of the paper.

\<{defn} \lbl{int-def} Let $N$ be  a structure, $M$ the union of some of the sorts of $N$.  $N$ is a {\em finite internal cover} of $M$ if 
$M$ is stably embedded in $N$, and $Aut(N/M)$ is finite (uniformly
in elementary extensions.)  Equivalently, $N \subseteq dcl(M,b)$
for some finite $b \in acl(M)$.  
 \>{defn}

A finite internal cover is a special case of an {\em internal cover}, where we demand that  $N/M$
is internal in place of $Aut(N/M)$ finite, and that $Aut(N/M)$ is definable.  
  In general, $Aut(N/M)$ is an $\infty$-definable group of $N$, isomorphic over $N$ to an 
 $\infty$-definable group of $M$.     cf. \cite{bedlewo}.

 While there is no difficulty in treating the general case,
 we will assume for simplicity of language that $Aut(N/M)$ is in fact {\em definable} in the
 internal covers considered in this paper.   (Only the case of   internal covers 
 with {\em finite} automorphism groups is needed for our applications.)

\<{remark} \lbl{lateruse} \rm Let $T''$ be a many-sorted expansion of $T$. For any $M'' \models T''$, let $M$ be the restriction to the language of $T$,
and let $Aut(M'') \to Aut(M)$ be the natural group homomorphism; let $K(M'',M)$ denote the kernel.  
  \<{enumerate}
\item If $Aut(M'') \to Aut(M)$ is always surjective, then $T$ is stably embedded in $T''$.  
\item If   in addition the kernel $K(M'',M)= Ker (Aut(M'') \to Aut(M))$ always has cardinality bounded in terms of $M$, then 
$T''$ is internal to $T$, i.e. all sorts of $T''$ are internal to the $T$-sorts. 
\item   If $K(M'',M)$ has cardinality bounded independently of $M$ and $M''$, then $T''$ is $T$-internal with 
an $\infty$-definable liaison group which is bounded, hence finite.  (Cf. \cite{bedlewo}, Appendix B.)  Thus in this case   each sort of $T''$ is a finitely imaginary sort of $T$. 
\item  If $Aut(M'') \to Aut(M)$ is always bijective, then $M'' \subset \dcl_{T''}(M)$.
\>{enumerate}
 \>{remark}

 The surjectivity implies that $T'$ induces no new structure
on the sorts of $T$, and also that $T$ is stably embedded in $T'$, cf. \cite{CH}, Appendix.  Injectivity of $Aut(M') \to Aut(M)$,  
implies that  $M' \subseteq \dcl(M)$; for suppose $c \in M'$ and $c \notin \dcl(M)$ .  We may take $M,M'$ to be sufficiently saturated
and homogeneous.   By  stable embeddedness there exists a small subset $A$ of $M$ such that $tp(c/A)$ implies $tp(c/M)$.  
As $c \notin \dcl(M)$, there exists $c' \neq c$ with $tp(c'/A)=tp(c/A)$.  So $tp(c/M)=tp(c'/M)$, and by stable embeddedness again
there exists $\si \in \Aut(M'/M)$ with $\si(c')=c$; the restriction of $\si$ to $M$ is the identity, but $\si$ is not, contradicting injectivity.  For the rest see  \cite{bedlewo}, Appendix B.

\<{lem} \lbl{no-expansion} Let $T'$ be a theory, $T$ the restriction of $T'$ to a subset of the sorts of $T'$,
$T''$ an expansion of $T'$ on the same sorts as $T'$.  Assume $T$ is stably embedded in $T'$,
and for any $N'' \models T''$, if $N',N$ are the restrictions to $T',T$ respectively,
the  natural map $Aut(N''/N) \to Aut(N'/N)$ is surjective.  Then $T',T''$ have the same
definable relations. \>{lem}

\prf  By the proof of Beth's implicit definability theorem, it suffices to show that $Aut(N'')=Aut(N')$ for any $N''$.  This is clear from the exact sequences 
$1 \to Aut(N'/N) \to Aut(N') \to Aut(N)$, $1 \to Aut(N''/N) \to Aut(N'') \to Aut(N)$, and
the equality $Aut(N'/N) = Aut(N''/N)$.  \eprf

 \>{section}
\<{section}{ Definable groupoids}
A {\em category} is a 2-sorted structure with sorts $O,M$, with maps $i_0,i_1: M \to O$
(the morphism $m \in M$ goes from $i_0(m)$ to $i_1(m)$), and a partial composition
$\circ:  M \times_{i_1,i_0} M \to M$, and an identity map $Id: O \to M$ (so that
$Id(x): x \to x$ is the identity map), satisfying the usual associative laws.  The language
of categories is thus 2-sorted, with relation  symbols $i_1,i_1,Id,\circ$.  

A (Grothendieck) {\em groupoid} is a category $\fG=(\ObG,\MorG)$ where every morphism has a 2-sided inverse.  For a groupoid $\fG$, let $Iso_{\fG}$ be the equivalence relation on $\ObG$:
$\MorG(c,c') \neq \emptyset$.  On the other hand, 
for any $a \in \ObG$, we have a group $G_a = \MorG(a,a)$.  These
groups are   isomorphic for $(a,b) \in Iso_{\fG}$:  : if $h  \in \MorG(a,b)$, then 
$x \mapsto h^{-1} x h $ is an isomorphism $G_b \to G_a$.  This isomorphism is well-defined up to conjugation.   Thus groupoids generalize, at different extremes, both groups and equivalence
relations:  
an equivalence relation is a groupoid with trivial groups, and a group is a groupoid with a single object.

We will 
assume   in this section that  $\fG$ has a unique
isomorphism type.  (I.e. 
$\ObG \neq \emptyset$, and $\MorG(a,b) \neq \emptyset$ for
all $a,b \in \ObG$. )  Without this assumption, one obtains relative
versions of the   results, fibered over the set of objects; for instance in \ref{g-triv}, the conclusion
becomes that one can interpret a set $S$ and a map $h : S \to T = \fG /\equiv$, such that for $t \in T$, for any representative $b \in \ObG$
of $t$, $F(b) $ is definably isomorphic to $S_t = h^{-1}(t)$.

If $X_a$ is a conjugation-invariant subset of some $G_a$, let $X_b= h^{-1} X_a h $,
where $h \in \MorG(b,a)$; the choice of $h$ does not matter.

In particular, if $N_a \triangleleft G_a$ is a normal subgroup, we obtain a system
of normal subgroups $N_b \triangleleft G_b$.  Moreover we can define
an equivalence
relation $N$ on  $\MorG (a,b) $:  $$(f,g) \in N \ \leftrightarrow \  g^{-1}f \in N_a \  
\leftrightarrow  \  fg^{-1} \in N_b$$ 
This gives rise to a   quotient groupoid with the same
set of objects, and with $\MorG ' (a,b) = \MorG (a,b) / N$.

It makes sense to speak of Abelian or solvable groupoids (meaning each $G_a$ is that.)

If $\ObG$ and $\MorG$ are defined by formulas in some structure $\Uu$,
as well as  the domain and range maps $\MorG \to \ObG$ and the
composition, we say that $\fG$ is a {\em definable groupoid} in $\Uu$.

A sub-groupoid is {\em full} if it consists of a subset of the objects,
with all morphisms between them.

Let $F: \fG \to Def(\Uu)$ be a functor.  We say that $F$ is {\em definable}
if $\{(a,d): a \in \ObG, d \in F(a) \}$ is definable, as well as 
$\{(a,b,c,d,e): a,b \in \ObG, c \in \MorG(a,b), d \in F(a), e \in F(b),
F(c)(d)=e \}$.

Similarly for $\star$-definable (= Pro-definable) or $\bigvee$-definable (see \secref{prelim}).
 But if there exist
a definable relation  $F_1$ and definable function $F_2$ such that for $a \in \ObG$, $F(a)=F_1(a)$,
and for $a,b \in \ObG$, $c \in \MorG(a,b)$, $F(c) = F_2(c)$, we will say
that $F$ is a (relatively) definable functor (even if if $\fG$ is only $\star$-definable.) 

\<{example}  \lbl{g-triv} Suppose each $G_a$ is trivial.  Then for each $a,b \in \ObG$
$\MorG(a,b)$ consists of a unique morphism.   In this case if $F:\fG \to Def(U)$ is a definable functor, one can interpret without parameters
a set $S$, definably isomorphic to each $F(a)$.    \rm 

Let
$$E_S = \{(a,b,a',b'): a,a' \in \ObG, b \in F(a), b' \in F(a'), 
\, \exists c \in \MorG(a,a'). \, F(c)(b)=b' \}$$
 $$S = \{(a,b): a \in \ObG, b \in F(a) \} / E_S$$

\>{example}

\<{example}   \lbl{ab}  If $\fG$ is Abelian, then the $G_a$ are all canonically
isomorphic, and one can interpret without parameters a single group, 
isomorphic to all $G_a$.    \>{example}
\proof  As in \ref{g-triv}:    the maps $G_a \to G_b$, being
unique up to conjugacy, are in this case in fact unique. \qed


\tsec{From $\star$-definable to definable groupoids}
\<{lem} \lbl{g-ext}  Let $\fG ^ 0$ be a  groupoid, with a distinguished
element $* \in \ObG ^0$.   Suppose $G^0_*=   \MorG ^0(*,*) $ is a subgroup of
a   group $G$.  Then $\fG ^0$ extends canonically to a groupoid
$\fG$ with the same objects, and with $  \MorG(*,*)=G$. \rm  

In other words, the natural map  
$\fG \mapsto \MorG(*,*)$,  from 
supergroupoids $\fG$ of $\fG ^0$ with the same object set, to supergroups $G$ of $G^0_*$, is surjective.   If $\fG^ 0, G$ are
$\star$-definable, so is $\fG$.  \rm  
 
\>{lem}

 { {\noindent \it Construction:} \hbox{       } }  Let 
 $$ (\MorG(a,b) = \MorG^0(*,b) \times G \times \MorG^0(a,*)) / \sim$$
where $(f,g,h) \sim (f',g',h')$ iff $((f') ^{-1} f) g   =g' (h'h^{-1})$.  Note that the
expression makes sense, since $((f') ^{-1} f), (h'h^{-1}) \in \MorG^0(*,*) \leq G$.
It defines an equivalence relation:  for instance, transitivity:  if
  $((f') ^{-1} f) g   =g' (h'h^{-1})$ and $ ((f'') ^{-1} f') g'   =g'' (h''(h')^{-1}))$,
then 
$$((f'')^{-1} f )g = ((f'') ^{-1} f')((f') ^{-1} f) g = ((f'') ^{-1} f')g' (h'h^{-1}) =
g'' (h''(h')^{-1})) (h'h^{-1}) = g'' (h'' h^{-1})$$  
Define composition 
$$\MorG(b,c) \times \MorG(a,b) \to \MorG(a,c)$$ by:
$ (j',g',f') / \sim) \circ  ((f,g,h)/ \sim) = (j', g'( f' f) g,h)/\sim$.

  The verifications are left to
the reader. \qed

\<{lem} \lbl{def-ext}  Let $\fG^ 0$ be a $\star$-definable  groupoid, with $G^0_a$
definable.  Then $\fG ^0$ extends to a definable groupoid $\fG$,
with $\MorG(a,b) = \MorG ^0(a,b)$ for $a,b \in \ObG^0$.

If $\sim$ an $\bigvee$-definable
equivalence relation, and all $a,b \in \ObG^0$ are $\sim$-equivalent, we can
obtain the same for $\fG$.   
\>{lem}

\proof  The hypothesis is intended to read:  $G^0_a$ is definable uniformly in $a$
(or equivalently, that the statement is true in any model.)  It follows
that $\MorG^0(a,b)$ is definable, for any $a,b \in \ObG^0$.  (This set
is a torsor over $G^0_a$, so it is definable with parameters; being 
$\star$-definable with parameters $a,b$, it must be definable uniformly in
these parameters.)  The definition of $\MorG^0(a,b)$ must extend over
all $a,b$ in some definable set $S_0$ containing $\ObG^0$.  The 
groupoid properties are certain universal axioms 
holding for all $a,b,c \in \ObG^0$; by compactness they must hold for all
$a,b,c \in S_1$ (some definable $S_1$, with $ \ObG^0 \subset S_1 \subset S_0$.)
Let $\ObG = S_1$, and use the definable function above to define $\fG$.

The two additional statements are also immediate consequences of  compactness. \qed

Some theories, notably stable ones (cf. \cite{meta}),   theories of finite S1 rank (\cite{pac}), and more generally supersimple theories (\cite{wagner}), have the property
that every $\star$-definable group is a projective limit $\lim_{\leftarrow n} G_n$,  where $G_1 \leftarrow G_2 \leftarrow \cdots$ is a sequence of definable groups and maps.  
As soon as this holds for one of the groups $G_a$, we can use the two lemmas above to
pass from an $\star$-definable groupoid to a definable one.


\tsec{The liaison groupoid}  \lbl{liaison}
 Let $\Uu$ be a universal domain for
 a theory admitting  elimination of quantifiers and
elimination of imaginaries.  

  Let $\Vv$ be a union of sorts of $\Uu$, closed under images of definable maps.  So $\Vv$ also admits elimination of imaginaries
.
 We will obtain a $\star$-definable groupoid; the set of objects
will be $\star$-definable, and the image of an $\infty$-definable set of $\Uu$; the sets of morphisms $\infty$-definable.  In
the situation of \ref{def-ext}, this will be the limit of definable groupoids.

\<{prop}  \lbl{groupoid}  Assume $\Vv$
is stably embedded in $\Uu$.

Let $Q$ be a   definable set of $\Uu$, internal to  
$\Vv$. 

   There exists   $\star$-definable groupoids $\fG$ in $\Uu$ and
$\fGV$ in $\Vv$, and   definable functors $F: \fG \to Def(\Uu)$ and
$F_{\Vv}: \fG_{\Vv} \to Def(\Vv)$, such that $\fGV$ is a full subgroupoid of $\fG$,
$F_{\Vv} = F | \fGV$, $\ObG = \ObGV \union \{*\}$, and $F(*)=Q$.

 We have  $F(G_*) = Aut(Q/\Vv)$.
  \>{prop}

\<{proof}  By internality of $Q$, and using   elimination of imaginaries in $\Vv$, there exists $Q_b$ definable over $b$ in $\Vv$,
and a $\Uu$-definable bijection $f_c: Q \to Q_b$.  Since $\Vv$ is stably 
embedded, $tp(c/b') \vdash tp(c/\Vv)$ for some $b'$; increasing $b$, 
we may assume $b=b'$.  (It is here that we must allow $b'$ to be
a tuple with an infinite index set.)

Let $\ObGV$ be the set of solutions of $tp(b)$, and let $\ObG =  \ObGV \union \{*\}$ (a formal element.)   Let
$\MorG(*,b') = \{c': tp(bc)=tp(b'c') \}$ (the morphism is viewed as identical
with the map $f_{c'}: Q \to Q_{b'}$.)  Let  $\MorG(b',*)$ be the same set
of codes, but each code viewed now as coding the inverse map
$Q_{b'} \to Q$.    Let $\MorG(*,*) = Aut(Q/\Vv)$.  

Observe the coherence of what has been defined so far: if $c,d \in
\MorG(*,b')$, then $tp(c/b') = tp(d/b')$.  Thus $tp(c/\Vv) = tp(d/\Vv)$.
Since $\Vv$ is stably embedded, there exists $\si \in Aut(\Uu/\Vv)$ with $\si(c)=d $.  Let $g = \si | Q$.  If $a \in Q$, then $\si(f_c(a)) = f_c(a)$
(since $\si$ fixes $\Vv$) but also $\si(f_c(a)) = f_{\si(c)}(\si(a)) = f_d(g(a))$.  Thus $f_c = f_d \circ g$.  Conversely, if $\tau \in Aut(\Uu/\Vv)$
is arbitrary, $h = \tau | Q$, then $f_d \circ h = f_{\si^{-1}(d)}$.

For $b',b'' \in \ObGV$, let $\MorG(b',b'')$ be the set of maps $Q_{b'} \to Q_{b''}$ of
the form $f_{c''} \circ {f_{c'}}^{-1}$, where $c' \in \MorG(*,b')$,
$c'' \in \MorG(*,b'')$.  

Note that if $\si(c'')=d''$, $g= \si | Q$, then  $f_{c''} = f_{d''} \circ g$,
$f_{c'} = f_{d'} \circ g$ for some $d' $ (=$\si(c')$), so
$f_{c''} \circ {f_{c'}}^{-1} =    f_{d''} \circ     {f_{d'}}^{-1} $.
Since $Aut(Q/\Vv)$ is transitive on $\MorG(*,b'')$,  fixing some
$d'' \in \MorG(*,b'')$, an arbitrary element of $\MorG(b',b'')$
can be written $ f_{d''} \circ     {f_{d'}}^{-1} $.  Similarly, 
an arbitrary element of $\MorG(b'',b''')$
can be written $f_{d'''}  \circ     {f_{d''}}^{-1} $.  So the composition 
of an arbitrary element of $\MorG(b',b'')$ with one of 
$\MorG(b'',b''')$ is an element of $\MorG(b',b''')$.  It follows that
we have indeed a groupoid.                   

Also, by expressing $\MorG(*,*) $ as $\MorG(b,*) \circ \MorG(*,b)$
for some $b$, it follows that $\MorG(*,*)$ is   an
$\infty$-definable set of permutations of $Q$ 
(over $b$, but a posteriori   over $\emptyset$, since at all events it is
invariant.)

Define the functor $F$ by $F(*)=Q$,
$F(b') = Q_{b'}$, and define $F$ on morphisms tautologically. 

Let $\fGV$ be the restriction of $\fG$ to $\ObGV$, and $F_V = F | 
\fGV$.   All the properties are then clear.
 \>{proof}

\<{remark} \lbl{groupoid-p} \begin{enumerate}

  \item  There exist definable maps $f_i: Q \to \Vv$ ($i \in I$)   such that $F(G_*)$ is transitive on each fiber of $f=(f_i)_{i \in I}$.

 \item  Assume   $G_* = \meet_n G_n$ where $G_1 \geq G_2 \geq \cdots$ are definable groups.  Then one can find 
a {\em definable} groupoid $\fG $ and a finite $I$  satisfying  Theorem \ref{groupoid} {\em except } the last statement and   (1) above.
 (And we still have $F(G_*)   \supset Aut(Q/\Vv)$).

\end{enumerate}
\>{remark}
\proof

\begin{enumerate}
  \item Let $f=(f_i)_{i \in I}$ enumerate all definable maps $Q \to \Vv$.  
  Then $c,d \in Q$ are $Aut(Q/\Vv)$-conjugate iff $tp(c/\Vv)=tp(d/\Vv)$ (by
  stable embeddedness) iff $f(c)=f(d)$.  
 
  \item  By \ref{g-ext}, \ref{def-ext}.

\end{enumerate}

\>{section} 

\<{section}{Generalized imaginaries}  \lbl{genim}

The  notion of an imaginary sort for a theory $T$ can be described as follows.
Let $T'$ be an extension of $T$ in a language containing the language
of $T$, and having an additional sort  $S$.   A universal
domain $\Uu'$ for $T'$ thus has the form $\Uu,S(\Uu')$.
$S$ is an {\em imaginary sort of $T$} if every model $M \models
T$ expands to a model $M' = (M,S_{M'})$ of $T'$ with
$S_{M'} \subseteq \dcl(M)$; equivalently (as noted above), for any such $M'$,
$$Aut(M') \to Aut(M)$$
is a group isomorphism. 

We will now consider    a slight generalization.       A  {\em finite generalized imaginary} sort
is defined as above, except that the   homomorphism
$$Aut(M') \to Aut(M)$$
is allowed to have finite kernel.  It is still assumed to be surjective.  More generally,
$S$ is called an {\em internal generalized imaginary sort} if the language of $T'$ 
 is  finite relative to the language of $T$ (i.e finitely many relation symbols are added),  and $T'$ is internal to $T$.     In this 
case, $Aut(S_{M'}/M)$ is isomorphic to $G(M)$ for some definable group $G$.    It makes sense to consider generalized sorts $S$ relative to a
sort $D$ of $T$, meaning that a definable map $S \to \bar{S}$ is given, and each fiber is an internal imaginary sort.  But  in this paper
we will consider internal generalized imaginary sorts almost exclusively, and will omit the adjective "internal".

An equivalent, more concrete definition of (ordinary) imaginaries can be given in terms of equivalence relations (cf. \cite{shelah}).  Let $E$ be a definable equivalence relation 
on a set $S$; then $S/E$ is added as a new sort, together with the canonical map $S \to S/E$.

This is used to find canonical parameters for definable families.  For $s \in S$, let $\delta(s)$
be a definable set; such that $\delta(s)=\delta(s')$ iff $(s,s') \in E$.  Then the image of $s$ in $s/E$
serves as a canonical parameter for $\delta(s)$.  

More generally, in place of equality,
one often has a definable bijection $f_{s,s'} : \delta(s) \to \delta(s')$, forming a commuting system.  Then for $\bar{s} \in S/E$ one introduces
$\delta(\bar{s} )$ as the quotient of the $(\delta(s): s/E = \bar{s})$ by the system $f_{s,s'}$,
  obtaining a canonical family $\delta(\bar{s}): \bar{s} \in S/E)$.   This can still be treated
  using equivalence-relation imaginaries, by an appropriate equivalence relation on $\dU_s \delta(s)$. However, if the system has more than one definable bijection $\delta(s) \to \delta(s')$, this fails.  We now generalize the above construction 
to more general groupoids.

 A {\em concrete definable category} of $T$ 
is a triple  $\fG = (\ObG, \MorG,\delta_{\fG})$ with $(\ObG,\MorG)$  a category
 interpretable in 
  $T$, and $\delta_{\fG} : \ObG \to Def(\Uu)$ a  faithful definable functor.

An {\em embedding} of $\fG_1$ into $\fG_2$ is a 0-definable fully faithful functor
$h: \ObG_1 \to \ObG_2$, together with a 0-definable system of definable
bijections $h_c: \delta_1(c) \to \delta_2( h(c))$ for $c \in \ObG _1$, 
such that $h_{c'} \circ  \MorG_1(c,c') = \MorG_2( h(c),h(c')) \circ h_{c}$.   In particular
conjugation by $h_c$ induces a group isomorphism $\MorG_1 (c,c) \to \MorG_2(h(c),h(c))$.

A {\em concrete groupoid} is a concrete category that is a groupoid.  

A   groupoid $\fG$  is {\em canonical} if $Iso_{\fG}$ is the identity,
i.e. two isomorphic objects of $\fG$ are equal.
A (concerete) groupoid $\fG$ is {\em a group (action)} if $\ObG$ has a single element.
 
Two concrete groupoids $\fG_1,\fG_2$ of $T$ are {\em equivalent} if there
exist $0$-definable embeddings $h_i: \fG_i \to \fG$ for some concrete groupoid
$\fG$, such that the image of $h_i$ meets every isomorphism class of $\fG$
(thus $h_i$ is an equivalence of categories.)
  In this case, $\ObG$ may be taken to be $\ObG_1 \du \ObG_2$, and 
the embeddings may be taken to be the identity maps.  If  $\fG_1,\fG_2$ and $\fG_2,\fG_3$
are equivalent, via concrete groupoid structures on $\fG_1 \du \fG_2$ and on $\fG_2 \du \fG_3$,
one may take the concrete groupoid generated by the union of these (with objects
$\ObG_1 \du \ObG_2 \du \ObG_3$) to see that $\fG_1,\fG_3$ are equivalent.

\tsec{The cover associated to a definable groupoid}
We describe a canonical cover of a theory $T$, associated with a definable groupoid $\fG$.  The new theory adds a distinguished object to
each isomorphism class of $\fG$.  
The cover will be internal if the groupoid
has a single isomorphism class.   A general groupoid $\fG$ can be viewed as a disjoint
union over $\nu \in \Ob \fG / Iso_{\fG}$ of the full sub-groupoid $\fG(\nu)$ whose objects are those of the isomorphism class $\nu$.   The cover $T'_{\fG}$ will then be   the 
  the free union of  the covers $T'_{\fG(\nu)}$.     The construction extends to the case of concrete groupoids.

Let $T$ be a theory, $\fG$ a definable groupoid, and $\dl: \fG \to Def(T)$ a definable
functor.  We construct a theory $T'=T'_{\fG,\dl}$ extending $T$.
 The   sorts of $T'$ are those of $T$, along with three new sorts $O,M,D$.  
The language of $T'$ is the language of $T$ expanded
by relations $i_0',i_1',\circ ', Id'$ for the language of categories on $O,M$ (with $O$ the objects, $M$ the morphisms),   and maps   $r: D \to O$,  $d: D \times_{r,i_0} M \to D$.    We will make $(O,M)$
into a concrete groupoid $\fG'$ with functor $\dl '$ by letting 
 $\dl'(x) = r \inv (x)$, and  $\dl'(f)$ be the restriction of $d$ to $r \inv (x) \times \{f\}$.   Finally the language
 has function symbols $j$ for a functor $(\fG,\dl) \to (O,M,\dl')$ of concrete categories.
 
The axioms of $T'$ are those of $T$, together with the statement that $(O,M,\dl')$ 
is a concrete groupoid,   $j$ is an embedding of concrete categories $\fG \to (O,M,\dl')$;
and $O$ has a unique element outside the image of $j$, in each isomorphism class.  

As usual we will write $Hom_{\fG'}(a,b)$ for $\{m \in : i_0(m)=a, i_1(m)=b\}$.  In particular $Hom_{\fG'}(a,a)$ forms a group, denoted
$Aut_{\fG'}(a)$.

\<{lem}  \lbl{cover} $T'$ is complete (relative to $T$).  $T'$  induced no new structure on the sorts of $T$.  
Each model $M$ of $T$ extends to a unique model $M'$ of $T'$, up to isomorphism
over $M$.    For any $a \in O(M')$, $\dl'(a)$ is internal to the sorts of $T$, and $Aut(\dl'(a) /T) = Aut_{\fG'}(a)$.  
 \>{lem}

\prf  Given $M \models T$, choose a representative $r_\nu$ of each isomorphism class
$\nu$ of $\Ob _\fG (M)$.   Let $O_0$ be  a copy  of $\Ob_{\fG}$, with $j: \Ob_{\fG} \to O$
a bijection; 
$*_\nu$ be a new element, $O = O_0 \union \{*_\nu: \nu\}$, and define a groupoid
structure in such a way that $j$ is an isomorphism of categories from $\fG$
the the sub-groupoid with objects $O_0$, and each $*_\nu$ is isomorphic
to each element of $\nu$.  It is easy to see that this can be done, and uniquely so
up to $M$-isomorphism.  In effect to construct $O$ one adds to each isomorphism 
class a new copy $*_\nu$ of $r_\nu$, and let $Mor(*_\nu, y)  $ be a copy of 
$Mor(r_\nu,y)$ for any $y \in \Ob \fG$, and $Mor(*_\nu,*_\eta) = Mor(r_\nu,r_\eta)$.
Similarly $\dl'(*_\nu)$ is a copy of $\dl(r_\nu)$.  For uniqueness, given two versions $O,O'$,
for any $\nu$ pick an isomorphism $f_\nu \in Mor(r_\nu,*_\nu)$, $f'_\nu \in Mor(r_\nu, *'_\nu)$,
and conjugate using $f$ from $*_\nu$ to $r_\nu$, then using $f'$ from $r_\nu$
to $*'_\nu$, to obtain isomorphisms $Iso(*_\nu) \to Iso(*'_\nu)$; compose
$\dl'(f') $ with $\dl' (f) \inv$ to obtain maps $\dl'(*_\nu) \to \dl'(*'_\nu)$; etc.

Completeness of $T'$ follows from the uniqueness of $M'$.  

Any element of $Aut_{\fG'}(a)$ acts on $\dl'(a)$, and also acts on any nonempty $Mor_{\fG'}(a,b)$ by conjugation; these combine to give a concrete groupoid automorphism fixing the image of
$j$, hence an automorphism fixing the $T$-sorts.  Given any automorphism
$\si \in Aut(T'/T)$, let $a \in M' \setminus M$ and pick $b \in \Ob \fG (M)$
with $a,b$ isomorphic in $\fG$.   we have $\si(a) =a$ since $\si(b)=b$ and
$a$ is the unique element outside the image of $j$ and isomorphic to $j(b)$.  
Pick an isomorphism $r \in Mor(a,b)$,  Then $\si(r)^{-1} r $ is an $\fG'$-isomorphism
of $a$, $\si$ coincides on $\dl' (a)$ and on any $Mor(a,c)$ with
the action of and conjugation by this element.

\eprf

{\bf Remark}  
The cover constructed above is 1-analyzable, i.e. relatively internal over a set interpretable in $T$, namely the set $Iso \fG$ of isomorphism classes of $\fG$; moreover
and has no relations among the fibers over $Iso \fG$.   In general, a 1-analyzable cover   $f:C \to D$ may have relations among fibers of $f$, not sensed by
the associated groupoid.  However any relation concerns finitely many fibers, so between them the groupoids associated to the   induced covers $f^n: C^n \to D^n$ for each $n$
do capture the information, and the cover may by coded by a definable simplicial groupoid.

{\bf Remark}  
 Assume $\fG$ has a single isomorphism class.  
if one fixes a parameter $b \in \ObG$, one may interpret the new element of $O$   by doubling $b$.  (Add one new object $b'$, and let
$Mor(b',c)$ be a copy of $Mor(b,c)$, etc., with the obvious rules.)
  In this
case, the corresponding groupoid imaginary is interpretable with parameters.  
However, unlike the groupoid imaginary sort,  this interpretation is incompatible with
 the automorphism group of the original structure.

%

\tsec{Internal covers and concrete groupoids}  
Two generalized imaginary sorts $S',S''$ of $T$ (with theories $T',T''$)   are {\em equivalent} if they are bijectively bi-interpretable over $T$, i.e. whenever $N' \models T', N'' \models T''$
are two models with the same restriction $M$ to the $T$-sorts, 
there exists a bijection $f: S_N \to S_{N'}$ such that
$f \union Id_M$ preserves the class of $0$-definable relations.

 As Levon Haykazyan and Rahim Moosa pointed out, 
 the theorem below requires an additional assumption, of {\em finite faithfulness}.    See their paper
 {\em  Functoriality and uniformity in Hrushovski's groupoid-cover correspondence,
Annals of Pure and Applied Logic, volume 169 (2018), number 8, 715--730. }
  
Finite faithfulness holds 
 automatically when  theory is stable
or when the automorphism groups in the groupoid are finite (or linear), as will always be the case in applications below.

\<{thm}  \lbl{im-gr}   
There is a bijective correspondence between internal imaginary sorts of $T$
and definable concrete groupoids   with a single isomorphism class (both up to equivalence.)
 \>{thm}
 
\prf  Given
the concrete definable groupoid $\fG$ with functor $F$, let $T'_{\fG}$ be the theory
described in \lemref{cover}.  Since $\fG$ has a single isomorphism class, 
there is a single element $*$ of $O$ outside the image of $\Ob_{\fG}$.    The sort
$S_{\fG}$ is taken to be $\dl'(*)$, with the structure induced from $T'$.
(Note that the rest of $T'$ is definable over the sorts of $T$ and $S_{\fG}$, using stable 
embeddedness.)

Conversely,  given an internal cover $N$, we obtain a *-definable concrete groupoid by \propref{groupoid}.  (The liaison groupoid of $N$.)
Since the number of sorts and generating relations is finite, it is clear that
$\Aut(S_N/M)$ is definable rather than $*$-definable.  By \lemref{def-ext} we can take
the groupoid $\fG$ definable.  It is   clearly
a concrete groupoid of $M$, well-defined up to equivalence. 

By \lemref{cover},  the liaison groupoid of $S_{\fG}$ is (equivalent to) $\fG$.  Conversely,
if we begin with $N$ and let $\fG$ be the liaison groupoid of $N$, then $S_{\fG}$
can be identified with $N$, though a priori $N$ may have more relations; but by construction
$S_{\fG}, N$ have the same automorphism group over $N$, so by \lemref{no-expansion}
their definable sets coinicide.    \eprf

In particular, a finite internal cover of $N$ may be realized as a  
  generalized imaginary sort, where the groupoid
  has a single isomorphism class, and finite isomorphism group at each point.

\<{example} \lbl{gext}  \rm Let $M$ be a finite structure, 
i.e.   finitely generated, with finitely many elements of each sort.
  Then any finite   extension 
$1 \to K \to \tG \to G \to 1$  of $G=\Aut(M)$
is the automorphism group of some finite  internal extension $\tM$ of $M$.  
The same holds in the $\aleph_0$-categorical setting, if the topology on $G$ is taken into account; see
\cite{ah-z} for proofs and \cite{ah-z1} for good examples.

\>{example}.

\<{lem}   Let $N$ be a finite internal cover of $M$, whose corresponding concrete groupoid
is equivalent to a (0-definable) group action.  Then the sequence
$$1 \to Aut(N/M) \to Aut(N) \to Aut(M) \to 1$$
is split.   \>{lem}

\proof In this case, the construction beginning with $\fG$ yields a structure
interpretable in $M$:  if $\Ob \fG = \{1\}$, the new structure has new sorts
$\delta(*)$ and $\delta(\Mor(*,1))$; by choosing a point of $\delta(\Mor(*,1))$ one 
obtains $M$ together with a copy $\delta(*) $ of $\delta(1)$ and a copy
$\delta(\Mor(*,1)) $ of $\delta(\Mor(1,1))$.  This interpretation yields a group homomorphism
$Aut(M) \to Aut(N)$ splitting the sequence.

This provides examples of structures that do not  eliminate  groupoid imaginaries.

\<{lem} \lbl{converse}  
Conversely, let $N$ be a finite internal extension of $M$, with associated  concrete groupoid $\fG$; and suppose the exact
sequence of automorphism groups is split functorially.  Then  $\fG$ is equivalent to a group action.  \>{lem}

\proof  By assumption, 
  there exists a subgroup $H \leq Aut(N)$, varying functorially when $(M,N)$ is replaced by an elementary extension, 
such that $H \to Aut(M)$ is an isomorphism.  Let $\bN$ be the expansion of $N$ by 
all $H$-invariant relations.   Then $\bN$ is bi-interpretable with $M$.  It follows
that the concrete groupoid corresponding to $\bN$ is equivalent to a group action.

\<{defn}
A finite internal cover $T'$ of $T$ is {\em almost split} if whenver $N' \models T'$,
with $N$ the restriction to the sorts of $T$, 
for some finite $0$-definable set $C$ of imaginaries of $N'$,
$N' \subseteq dcl(N \union C)$.  If $Aut(N'/C) \to Aut(N)$ is surjective,
we say that the cover is {\em split}.  
\>{defn}

Thus ``$T'$ almost split over $T$" is the same as: ``$T'_{\acl(\emptyset)}$ is split, over $\acl(\emptyset)$.''

\<{defn}  $T$ eliminates (finite, strict) generalized imaginaries if
every   concrete groupoid (with finite automorphism groups, with one isomorphism class)
 $\fG$ is equivalent to a canonical one.    \>{defn} 
 
Note that ordinary elimination of imaginaries holds iff every groupoid with trivial groups
is equivalent to a canonical one.

\<{lem}  $T$ eliminates  finite generalized imaginaries iff $T$ eliminates finite imaginaries, and
 every finite internal cover  of $T$ is split.   

  \>{lem} 

\prf  We use \thmref{im-gr}.  $T$ eliminates   finite generalized imaginaries iff 
every concrete groupoid $\fG$ with finite automorphism groups and one isomorphism class
is equivalent to a group action.  If $\fG$ is a group action, 
 the finite internal cover corresponding to $\fG$ is clearly split.  Conversely if 
 the cover $T'$ of $T$ is split, it has an expansion $T''$ bi-interpretable with $T$.
 $T''$ is still a finite internal cover, and by \thmref{im-gr} corresponds to a sub-groupoid $\fG''$
 of $\fG$, with one isomorphism class and trivial automorphism groups.  Let $*$ be a formal element corresponding
 to the isomorphism class of $\fG$.  We may assume $\Ob \fG = \Ob \fG'$.
 For $a \in \Ob \fG$, let $Mor(a,*) = Mor_{\fG}(a,a)$.
 Given $a,b \in \Ob_{\fG}$, there is a unique $f_{a,b} \in \Mor_{\fG'}(a,b)$.  Use 
 $\dl(f_{a,b})$ to identify $\dl (a),\dl (b)$, and let $\dl (*)$ be the quotient.  Also use $f_{a,b}$
 to identify $\Mor(a,a)$ and $\Mor(b,b)$, by composition, and let $\Mor(*,*)$
 be the quotient.  We have found a common extension of the group action 
 of $\Mor(*,*)$ on $\dl (*)$, and of the concrete groupoid $(\fG,\dl )$.  \eprf



\<{rem} \lbl{modelcov} If algebraic points form an elementary submodel of $M$,
then every finite internal cover is almost split.   \rm  Indeed by definition, a finite internal cover $N$ satisfies $N \subset \dcl(M,C)$ with 
$C$ a finite, $M$-definable set.  As $\acl(\emptyset) \prec M$, we can choose $C \subset \acl(\emptyset)$.\>{rem}

A definable group homomorphism $f: \tH \to H$ is a {\em definable central extension}
if $f$ is surjective and $\ker{f}$ is contained in the center of $H$.  We now relate
finite internal covers of internal covers of a theory $\tau$ to definable   central extensions 
of the liaison  group of the latter.   Assumption (3) below says that finite generalized
imaginaries of $M$ arising from definable finite central extensions of groups 
are eliminable; the conclusion is that all finite generalized imaginaries are.

\<{prop}  \lbl{finitecentral} Let $T$ be a a theory with   a distinguished stably embedded
sort $\k$, $\tau=Th(\k)$.     Let $M \models T$.  Assume $T$ eliminates imaginaries, and:

\begin{enumerate}
  \item  Every finite internal cover of $\k$  is almost split.
  \item  Let $D$  be  a $T$-definable  set.  Then $D$ is $\k$-internal, and
      $Aut(D/\k)$ is $T_M$-definably isomorphic to a  $\tau$-definable group $H$. 
  
  \item Let $F_0$ be a finite definable set of imaginaries of $T$, 
  $D$ a $T_{F_0}$-definable set, $h: Aut(D/ \k,F_0) \to H$ an $M_T$-definable group
  isomorphism (as in (2)).  Let  $f: \tH \to H$ be a $\tau$-definable central extension.  Then 
   there exists  a finite $T$-definable $F$ containing $F_0$, and a
 $T_F$- definable $\tD $ containing $D$,  and injective $T_M$-definable group homomorphisms
 $\th: Aut(\tD/k,F) \to \tH$, $h: Aut(D/k,F) \to H$, with images of finite index, and
 $f \th = h$.
 
 \item For any $\tau$-definable group $H$, and any finite Abelian group $K$, the group of $M$-definable homomorphisms $H \to K$ is finite.
   \end{enumerate} 

 Then any finite internal cover of $M$ is almost split.  \>{prop}

\<{proof}  Let $M' = M \union \{C\}$ be  a finite internal cover  of $M$, with $C$ definable;
$T'=Th(M')$.  
We have $G:=  Aut(M'/M) = Aut(C/M) = Aut(C/D)$ for some definable set $D$ of $T$; 
$G$ is  a finite group.  We may enlarge $C$ so that $D \subseteq C$. 

Two preliminary remarks:

If $F$ is a finite definable set of imaginaries of $T'$, there exist finite definable sets
$F_T$, $F_\tau$ of imaginaries of $M,\k$ respectively, such that for for any $M,M'$ as above,
$\dcl(F) \meet M^{eq} = \dcl(F_T) \meet M^{eq}$, $\dcl(F) \meet \k^{eq} = \dcl(F_\tau) \meet \k^{eq}$.  Thus $Aut(M'/F) \to Aut(M'/F_T) \to Aut(\k_M / F_\tau)$ are defined and surjective.
Thus to show that $T'$ is split, it suffices to prove the same for $T'_F$.   In particular,
taking $F$ to be the finite group $G$, we may assume each element of $G$ is 0-definable.
In this case, $G$ is central in $Aut(M'/\k)$.   We have a central extension
$$1 \to G  \to Aut(C/\k) \to Aut(D/\k) \to 1$$ 
By internality, the sequence is isomorphic
to a central extension $$1 \to K \to \tH \to H \to 1$$ 
of $\tau$-definable groups, via an $M'$-definable map  $\tf: Aut(C/\k) \to \tH$, and
an $M$-definable map 
$f: Aut(D/ \k) \to H$.  

The condition in (3) is stated for central extensions with prime cyclic kernel;   by iteration it is closed under all finite central extensions.  

Hence, after naming parameters for a further finite definable set, and passing
to corresponding subgroups of finite index in $\tH,H$, there exists
a $T$-definable $\tD$ (containing $D$) such that 
$Aut(\tD/ \k) \to Aut(D/ \k)$ is isomorphic to $\tH \to H$, by $T$-definable maps
$\th, h$; and $h = f$.  
 
 Now $\tH \times_H \tH$ has a subgroup of finite index isomorphic to $\tH$, namely
the diagonal subgroup $\Delta_{\tH}$.  $\Delta_{\tH}$ is invariant under any $M$-definable
automorphism of $\tH \times_H \tH$ of the form $(\alpha \times_\beta \alpha)$,
with $\alpha: \tH \to \tH $ an automorphism lying over $\beta: H \to H$.  
But any automorphism of $\tH$ over $H$ has the form $x \mapsto z(x) x$ for
some homomorphism $z: \tH \to K$.  By (4), there are only finitely many such definable
homomorphisms, and so the group of $M$-definable automorphisms of $\tH$ over $H$ 
is finite, and hence the automorphisms  $(\alpha \times_\beta \alpha)$ have
finite index within the group of all $M$-definable automorphisms of $\tH \to H$.  So
$\Delta_{\tH}$ has finitely many conjugates by such automorphisms.  Taking their
intersection, we find a subgroup $S$ of $\tH \times_H \tH$ of finite index, mapping
injectively to each factor $\tH$, and invariant under all $M$-definable
automorphisms of $(\tH,H)$.  

It follows that the pullback $S'$ of $S$ under $(\th,\tf,f)$ does not depend on the choice
of the triple $(\th,\tf,f)$.  It is thus a definable subgroup of   
$Aut(\tD / \k) \times_{Aut(D/\k)} Aut(C / \k)$.  Hence  $Aut(\tD,C / \k)$ also has a definable subgroup of finite
index $S''$ mapping injectively to $Aut(\tD/\k)$.  
 
By \cite{bedlewo}, there exists a $T'$-definable set $Q$ with $Aut(Q / \k) = Aut(\tD / \k)$,
and such that   $Aut(\tD / \k)$ acts transitively on $Q$,
with trivial point stabilizer.  The quotient $Q / S''$ is a  finite internal cover of $\k$.   By (1), for some $0$-definable
 finite set $F'$ we have $Q/ S'' \subseteq \dcl(F')$.  So $Aut(M' / F')$ is contained
 in $S''$.  Hence $Aut(M' / F', M) \subseteq Aut(M' / F',\tD) = 1$.  So $M'$ is an almost split extension of $M$.   
\>{proof}

 \tsec{Groupoids in ACF}
Consider $ACF_L$, the theory
of algebraically closed fields containing a field $L$.   Let $L^a$ be the algebraic closure of $L$.       Every 
concrete groupoid is equivalent to a a subgroupoid with finitely many objects
in each  equivalence class.  The question essentially reduces to concrete groupoids
with finitely many objects.  

If $L$ is real-closed, the Galois group $Aut(L^a/L)$ is $\Zz/2\Zz$, and admits nontrivial central extensions $1 \to Z \to E^a \to Aut(L^a/L) \to  1$.  $E^a$ lifts to an extension $E$ 
of $Aut(K/L)$ (where $K \models ACF_L$.)   The sequence
$1 \to Z \to E \to Aut(K) \to 1$ is not split, any more than the $E^a$-sequence.
 As in \exref{gext} there exists a finite internal cover  $M$ of $ACF_L$ with $Aut(M)=E$.
 The concrete groupoid corresponding
to $M$ cannot be  equivalent to a canonical one.

On the other hand, if $L$ is PAC, then (cf. \cite{FJ}) $Aut(L^a/L)$ is a projective profinite
group.  In this case every finite concrete groupoid should be equivalent to a canonical one.

\<{problem}  Give a geometric description of  the groupoid-imaginaries when 
when $L$ is a finitely generated extension of an algebraically closed field. \>{problem}

\>{section} 

\<{section}{Higher amalgamation}  

Let $T$ be a theory (or Robinson theory), for simplicity with quantifier elimination.  A $T$-structure is an algebraically closed substructure of a model of $T$.  
Let $\fC_T$ be the category of algebraically closed $T$-structures.
A partially ordered set $P$ can also be viewed as a category, and
we will consider functors $P \to \fC_T$.   Specifically let
$P(N)$ be the  partially ordered set of all subsets of $\{1,\ldots,N\}$,
and let $P(N)^-$ be the sub-poset of proper subsets.

By an {\em $N$-amalgamation problem} we will
mean a functor $A: P(N)^- \to \fC_T$.   A {\em solution} is a functor $\bar{A}: P(N) \to \fC_T$,
where $P(N)$ is the partially ordered set of all subsets of $P(N)$,
extending $A$.  We will 
 demand for both $a=A,\bar{A}$ that $a(s) = \acl \{ a(i): i \in s \}$ (this is by no means essential, but
 simplifies the definitions of independence-preservation and of uniqueness of solutions below.)   
 
We assume $T$ is given with a notion of {\em canonical $2$-amalgamation}.  I.e. we are given a functorial solution of all 2-amalgamation problems.  Equivalently, we have a notion of {\em independence} of two substructures of a model of $T$,  over a third; or again, 
a functorial extension process $p \mapsto p|B$ of types over $A$ to types over $B$,
 where $A \leq B \in \fC_T$.  
We  assume that this notion of independence is symmetric and transitive,
 cf. \cite{baldwin}.    When for any $B$, 
  $p|B$ is   an $A$- definable type, we will say that amalgamation is definable at $p$.
  (This is always the case for stable theories, cf. \cite{shelah}.) 
   
(The ``uniqueness of non-forking extensions'' comes with the presentation here; some of
the considerations below generalize  easily  to the case of a canonical {\em set} of solutions rather than one.)

A functor $A: P \to \fC_T$ is {\em (2)-independence-preserving} if it is 
compatible with the given canonical $2$-amalgamation; i.e. 
 whenever $s = s_1 \meet s_2 \subset s' \in P$,
$A(s_1),A(s_2)$ are independent over $A(s)$ within $A(s')$.

At this point, we consider the problem of {\em independent amalgamation}.   An independent amalgamation problem (or solution) is a functor $A: P \to \fC_T$ (where $P =  P(N)^-$,  respectively $P= P(N)$) 
compatible with the given canonical $2$-amalgamation; i.e. 
such that whenever $s = s_1 \meet s_2 \subset s' \in P$,
$A(s_1),A(s_2)$ are independent over $A(s)$ within $A(s')$.   We will also
demand:  $A(\emptyset) = \acl(\emptyset)$.

Let us say that $T$ has $n$-uniqueness (existence,exactness) if every independent $n$-amalgamation problem has at most one (at least one,
exactly one) solution, up to isomorphism.    

Similar diagrams appear  in work of Shelah in various contexts, cf. e.g. \cite{sh87b}.
Elimination of imaginaries was introduced in \cite{shelah} precisely in order to obtain
$2$-exactness for stable theories.     $3$-existence follows, but   $4$-existence,
and $3$-uniqueness, can fail:  cf. \cite{kpy}.  We will see below, however, 
that with generalized imaginaries taken into account, stable theories are $3$-exact.

Occasionally we will also require 
 $(n-1,n+k)$-existence for $k \geq 1$.  This means that a solution exists to every
{\em partial} independent amalgamation problem
 $(a(u): u \subset \{1,\ldots,n+k\}, |u| < n \}$.  We have however:
 
 \<{lem} \lbl{weakstepup} Assume $T_A$ has $n$-existence for all $A$.   Then 
 
 (1)  $T$ has   $(n-1,n+k)$-existence for any $k \geq 1$.  
  
 (2) If $n$-uniqueness holds, so does  $n+1$-existence.
  \>{lem}

\prf  (1) An easy induction.  For instance take $k=1$, and assume given  $(a(u): u \subset \{1,\ldots,n+1\}, |u| < n )$.   
Let $U = \{u \subset  \{1,\ldots,n+1\}: |u| \leq n, (n+1) \in u \}$.  For $u \in U$, 
by $(n-1,n)$-existence, one can find $a^*(u)$ extending $a(v)$ for $v \subset u, |v|<n$.
But $U$ 
is isomorphic
to the set of subsets of $\{1,\ldots,n\}$ of size $<n$, so by another use of 
$(n-1,n)$-existence $(a^*(u): u \in U)$  admits a solution $(b(u): u \subseteq \{1,\ldots,n+1\})$; this also solves the original problem $a$.

(2).  We will use $(n-1,n+1)$-amalgamation.  
    Given an $n+1$-independent amalgamation 
 problem $b$, let $a$ be the restriction to the faces $u$ with $|u|<n$.  This problem
 has a solution $c$; for each $u$ with $|u|=n$,  $c$ restricts to a solution $c(u)$ of the problem 
 $(a(v): v \subset u)$; by $n$-uniqueness, these solutions must be isomorphic to the
 original solutions $b(u)$.  By means of these isomorphisms $b(u) \to c(u)$ ($|u|=n$),
 $c$ provide a solution to the problem $b$.
  \eprf


\<{lem} \lbl{3ul}  Let $T$ be a theory with a canonical $2$-amalgamation, admitting 
elimination of imaginaries.  For any independence-preserving functor $a: P(3) \to \fC_T$, 
the following two conditions are equivalent:  

(1)    $a(12) \meet \dcl(a(13),a(23)) = dcl(a(1),a(2))        $
 
(2)   If  $c \in a(12) = \acl(a(1),a(2))$, then
$tp(c/a(1),a(2))$ implies $tp(c / a(13),a(23))$.

Moreover, 
$3$-uniqueness is equivalent to the truth of (1,2) for all such $a$.   

  \>{lem}

\proof  Assume (1) holds.  In the situation of (2),  the solution set $X$ of 
$tp(c/a(13),a(23))$ is a finite set, hence coded in $a(12)$, and
defined over $a(13) + a(23)$; thus by
(1), $X$ is defined over $a(1),a(2)$; being consistent
with $tp(c/a(1),a(2))$, it must coincide with it. 

Conversely, if $c \in a(12) \meet \dcl(a(13),a(23)) $, then $tp(c / a(13),a(23))$ has
the unique solution $c$, so if (2) holds then the same is true of $tp(c/a(1),a(2))$, 
and hence $c \in \dcl(a(1),a(2))$.

Suppose (1) fails.  Then the restriction map
$$Aut(a(12) / a(13),a(23)) \to Aut(a(12) / a(1),a(2))$$
is not surjective.  Let  $\si \in Aut(a(12) / a(1),a(2))$ be an automorphism
that does not extend to $a(13)a(23)$.  Let $a'$ be the same as $a$ on subsets
of $\{1,2,3\}$, and also the same on morphisms except for the inclusion 
$i: \{1,2\} \to \{1,2,3\}$; and let $a'(i) = a(i) \circ \si$.  Then $a'$ is a solution
to the independent amalgamation problem $a | P(3)^-$, and is not isomorphic to $a$.  So 3-uniqueness fails.  

Conversely, suppose we are given an  independent amalgamation problem 
$a : P(3)^- \to \fC_T$, and two solutions $a',a''$ on $P(3)$.  We may take $a'$ to take
the morphisms to inclusion maps.  Then all $a(ij)=a'(ij)$, and are embedded
in $A=a'(123)$.  We can identify $a''(123)$ with $A=\acl(a_1,a_2,a_3)$ also.
Then the additional data in $a''$ consists of isomorphisms $a(i,j) \to a(i,j)$,
compatible with the inclusions of the $a(i)$.   By 
$2$-uniqueness, we may further assume that these isomorphisms are the identity on
$a(2,3)$ and $a(1,2)$; so that $a''$ reduces to an automorphism $f$ of $a(1,3)$,
over $a(1),a(3)$.  By (1), $f$ extends to an elementary map fixing $a(12),a(13)$.  
This further extends to an automorphism $F$ of $a(1,2,3)$.  $F$ shows that
the two solutions $a',a''$ of the problem $a$ are isomorphic.   
 \qed

\<{prop} \lbl{cond}   Let $T$ be a stable theory admitting elimination of quantifiers and of imaginaries.   
Assume  every finite internal cover of $T_A$ 
    almost    splits over $A$.   Then $T$ has 3-uniqueness.
 \>{prop}
 
In place of stability, we can assume $T$ is given  with a notion of   2-amalgamation, 
and show \lemref{3ul} (1) holds whenever the amalgamation is definable at one of the vertices
of the triangle in question. 
 
We will see later that $4$-existence  is equivalent to 3-uniqueness.   

Compare  \cite{bouh}, where a finite internal cover  was constructed in the same way;
the purpose there  was to interpret a group from the group configuration, in a stable theory.
This is also done for simple theories in  \cite{kpy}, where 4-existence is assumed.
In hindsight, it all coheres.

\proof  
 Since by definition  $3$-uniqueness for $T_{\acl(\emptyset)}$ implies $3$-uniqueness for $T$,
 we may assume $\acl(\emptyset)=\dcl(\emptyset)$.
 
 Let $a: P(3) \to \fC_T$ be an independence-preserving functor, with notation as in
\lemref{3ul}.  
Replacing $T$ by $T_{a(\emptyset)}$ we may
assume $a(\emptyset) = \acl(\emptyset) = \dcl(\emptyset)$.  Fix an enumeration of $a(i)$.   We will describe
a finite internal cover   $T^+$ of $T$, associated with $a$.

Let $FU$ be the set of formulas $S(x_1,x_2;u)$ such that whenever $M \models S(a_1,a_2;c)$, 

\begin{enumerate}
  \item $c \in \acl(a_1,a_2)$ 
  \item If $tp(a_3) = tp(a(3))$ and $a_3$ is independent from $\acl(a_1,a_2)$, then
  $c \in \dcl(\acl(a_1,a_3), \acl(a_2,a_3))$. 
 \end{enumerate}

If $a_i$ enumerates $a(i)$, then 
by definability of the canonical extension of $tp(a(3))$, for any $c \in a(12) \meet \dcl(a(13),a(23))$
there exists $S \in FU$ with $S(a_1,a_2,c)$.  

Let $S,S',S'',\ldots  \in FU$.
By definability of the canonical extension of $tp(a(1))$, for any formula 
$\phi(x,y,z,y',z',y'',z'',\ldots,w)$ there exists  a formula $\phi^*(y,z,y',z',\ldots,w)$
(depending on $\phi$ and on the sequence $S,S',\ldots$)
such that for any $b,c,b',c',\ldots,d$, any any $a \models tp(a(1))$ with 
$a$ independent from $\{b,c,b',c',\ldots,d\}$ and such that $S(a,b,c),S'(a,b',c'),\ldots$,
$$\phi(a,b,c,b',c',\ldots,d) \iff \phi^*(b,c,b',c',\ldots,d)$$

We construct a many-sorted cover $T'$ of $T$ as follows.  

Let $L^+$ be a language containing $L$, as well as a new sort
$NS(u)$ for any $S \in FU$; and a definable map $f_S$; and for each $S,S',S'' \ldots \in FU$
and each $\phi(x,y,z,y',z',\ldots,w)$, a relation $N\phi(z,z',\ldots,w)$.   Given $M \models T$,
we construct an $L^+$ -structure $M^+$ as follows.  Within some elementary extension
$M^*$ of $M$, let $a \models  tp(a(1)) | M$.
Let $$NS(M^+) = \{(b,c): b \in M, M^* \models S(a,b,c) \}$$
and let $f_S(y,z) = y$.  For $d \in M, c \in NS(M^+), c' \in NS'(M^+), \ldots$,
interpret $N\phi$ so that
$$N\phi((b,c),(b',c'),\ldots,d) \iff \phi(a,b,c,b',c',\ldots,c)$$

Using the definability of $tp(a(1)) | M$, one sees that $T^+ = Th(M^+)$ does not depend on any of the choices made.

We now use \remref{lateruse}.    Each sort of $T^+$ will
be seen to be a finite internal cover of $T$, as soon as we show:

{\medskip \noindent \bf  Claim. \ }  $T^+$ is a bounded internal cover of $T$.

\proof  Given $M \models T$, we constructed an expansion $M^+ \models T^+$ of the
same cardinality, such that $Aut(M^+) \to Aut(M)$ is surjective.  It remains to show that the kernel is bounded.  $M^+$ can be constructed as follows.
We have $M \models T$.  Let $M_3$ be an elementary extension of $M$,
with $a_3 \in M_3$, $a_3 \models tp(a(3)) | M$.  Let $M^*$ be an elementary extension
of $M_3$, with $a_1 \in M^*$, $ a_1 \models tp(a(1)) | M_3$.  We can construct $M^+$
using $M,a_1$, so that $M^+ \subseteq \acl(M,a_1)$;
actually $M^+ \subseteq \union_{a_2 \in M} \acl(a_1,a_2) \meet \dcl(\acl(a_2,a_3),\acl(a_1,a_3))$.  
  Any automorphism of $M^+$ over $M$ lifts to an automorphism of $M_3^+$ over $M_3$,
 which in turn is an elementary automorphism of $M_3^+$ (viewed as a subset of $M^*$)
 over $M_3(a_1)$.    Thus the homomorphism 
$Aut(\dcl(acl(a_1,a_3),M_3 )/ M_3,a_1,a_3) \to Aut(M^+)$ is surjective.  But the first group is clearly bounded.    \qed{(Claim)}

Now by assumption, every finite internal cover of $T $ is almost split.
Let $M$ be a model of $T$ containing $a(2)$, and let $a(1)$ be independent from $M$,
in some elementary extension $M^*$ of $M$.  Then $M^+$ can be embedded into $\dcl( a(1),M)$.
Let $c \in a(12) \meet \dcl(a(13),a(23))$, so that $c \in M^+$.  
Then $c \in \dcl( a(1), M)$.  But $tp(a(1),c/M)$ is $a(2)$-definable, since 
$tp(a(1)/M)$ is definable, and $c \in \acl(a(1),a(2))$, and $a(2)$ is algebraically closed.
Thus $c \in \dcl(a(1),a(2))$.  This proves the property of \lemref{3ul} (1).
\qed

\ssec{Adding an automorphism}  
\lbl{add-automorphism}

  We include a general lemma on adding an automorphism to a stable theory,  that will aid in describing the linear  imaginaries of pseudo-finite fields.   This was the route taken in 
  \cite{pac} to the imaginaries of the pseudo-finite fields themselves;
  it appears best to repeat it from scratch in the linear context.  
In \cite{pac}, as here, only the fixed field was actually needed.  The imaginaries for the full theory were considered (and eliminated) in \cite{CP} for strongly minimal $T$ (in \cite{CH} for $T=ACFA$).
 An unpublished example of Chatzidakis and Pillay shows
that it is not true in general.   We show however that the principle is correct if
  generalized imaginaries are taken into account.

   Let $T$ be a   theory with  elimination of 
quantifiers and elimination of imaginaries.
(In our application, $T$ will be a linear extension of the theory of 
 algebraically closed fields.)

Let $\tC = \{(A,\si): A \in \fC_T, \si \in Aut(A) \}$.  We define independence for $\tC$
by ignoring the automorphism $\si$.   In the present framework, $2$-uniqueness will not hold; this is because of the choice
involved in extending an automorphism from $\dcl(A_1\union A_2)$ to $\acl(A_1 \union A_2)$.

Consider pairs $(A,\si)$, with $A$ an algebraically closed substructure
of a model of $T$, and $\si$ an automorphism of $A$.  This is the class of models of a theory $T^\forall _\si$, in a language where quantifiers over $T$-definable finite sets
are still viewed as quantifier-free.  
  Under certain conditions, including the application in \S4 to linear theories over ACF, 
$T^{\forall}_\si$ has a model completion, a theory $\tTs$ whose
models are the   existentially closed  models of $T^{\forall}_\si$.
$\tTs$ is unique if it exists.  At all events, $\tC$ amalgamates to a universal domain, and can be viewed as a Robinson theory.

\<{prop}\lbl{su}  Let $T$ be a theory with a canonical $2$-amalgamation, admitting 
elimination of imaginaries.  Assume $T_A$ has   $n$-existence
  Then 
conditions (1)-(4) are equivalent.  
 \begin{enumerate}
  \item $n$-uniqueness.

  \item $n$-existence for $\tC$
  \item  Let  $a: P(n)^- \to \fC_T$ be independence-preserving.  Let 
    $u_0 = \{1,\ldots,n-1\}$;  let $a(<v) = \dcl( \{a(v'): v' \subset v, v' \neq v\})$,
  $a(\not \geq v) = \dcl(\{a(v'): v   \not \subseteq v'\})$.
  then
  $$a(u_0) \meet a(\not \geq u_0) = a(< u_0)$$
  \item  With $a,u_0$ as in (3), 
   $$\Aut(a(u_0) / a(\not \geq u_0)) = \Aut(a(u_0) / a(< u_0))$$
 
\end{enumerate}
 \end{prop}
 
 \prf 
(1) $\implies$ (3) is proved as in \lemref{3ul}.  
 
(3) $\iff$ (4):   Using imaginary Galois theory, cf. \cite{poizat}).

(4)  $\implies$ (1):  Let $a: P(n)^- \to \fC_T$ be an $n$-amalgamation problem,
and let $a',a''$ be two solutions.  As in  \lemref{3ul} we may assume 
that for each $u \subset n$, $a'(Id_u)$ is the inclusion  of $a(u)$ in 
$\acl(a(1),\ldots,a(n))$.  Now for $i \in \{1,\ldots,n\}$, $u^i= \{1,\ldots,n\} \setminus \{i\}$,
$Id_{u^i}$ the inclusion of $u^i$ in $\{1,\ldots,n\}$, $a''(Id_{u^i})$ 
is an isomorphism $a''(u^i) \to a''(\{1,\ldots,n\})$, i.e. an automorphism $f^i$ 
of $a(u^i) = \acl( a(j): j \in u^i )$; and since $a''$ extends $a$, 
$f^i \in  \Aut(a(u^i) / a(< u^i))$.  By (4), 
$f^i \in \Aut(a(u^i) / a(\not \geq u^i))$.  So $f^i$ extends to an automorphism $F^i$ 
of $a(\{1,\ldots,n\})$ fixing $a(\not \geq u^i))$.  Let $F$ be the product of the $F^i$
(choose any ordering.)  Then $F | a(u^i) = f^i$.  So $F$ shows that $a'',a'$ are isomorphic.

$(1) \implies (2)$:   consider an independent $n$-amalgamation problem for $\tC$; it consists of an independent $n$-amalgamation
 problem $a = (a(u): u \subset \{1,\ldots,n\})$ and a compatible system of
 automorphisms $\si_u \in \Aut( a(u))$.  Using $n$-existence, extend $a$ to a solution;
 it is a system $(b(u): u \subseteq \{1,\ldots,n\})$, and compatible embeddings
 $f_u: a(u) \to b(u)$.  Now let $g_u = f_u \circ \si(u): a(u) \to b(u)$.  Then $(b,g)$
 is another solution.  By $n$-uniqueness, the two solutions must be isomorphic;
 so there exists $\si: b(\{1,\ldots,n\}) \to b\{1,\ldots,n\})$
 such that $g_u = \si f_u$.  This shows exactly that $(b,\si)$ is a solution to the 
 original automorphic problem, via $f$.  
 
$(2) \implies (3)$:  Let  $\si_0 \in \Aut(a(u_0) / a(< u_0))$,
let $\si_u = Id_{a(u)}$ for every $u \not \geq u_0$.  View this data as an independent
amalgamation problem for $\tC$.  By (3), it has a solution $a'$.  We use 2-uniqueness
to note that $tp(a'(u_0),a'(n)) = tp(a(u_0),a(n))$.  Thus $\si_0$ has an extension 
to $Aut(a(\{1,\ldots,n\}))$ fixing $a(\not \geq u_0)$.  So $\si_0$
fixes $a(u_0) \meet a(\not \geq u_0)$.  By imaginary Galois theory again,
$a(u_0) \meet a(\not \geq u_0) \subseteq a(< u_0)$.


\eprf



\<{remark}  If $n$-uniqueness fails, it fails already for the Abelian algebraic closure.  
For the Abelian algebraic closure, a formulation in terms of  homological algebra becomes possible.  \>{remark} Fix types $p_1,p_2,\ldots$ or more simply one type $p$.  Let $a_1,\ldots,a_k$
be $k$ independent realizations of $p$.  Let $G_k$ be the Abelianization of 
$Aut( \acl(a_1,\ldots,a_k) / \dcl(a_1,\ldots,a_k))$.  There is a natural homomorphism
$G_{k+1} \to G_k^{k+1}$, restricting to each $k$-face.  In terms of this basic data, one can 
describe  homologically the questions of $n$-existence and uniqueness.   The point is
that in the independence preserving functors $a$, the   $a(u)$ can be taken
to be standard objects $\acl(a_i: i \in u)$, so only the image of morphisms under the functor
matters.


The proof of
\propref{stableEI} below follows the same outline as \cite{pac},\cite{CH}, \cite{CP}.  It may be possible to give
a proof based on minimizing $tp(a/A)$ in the fundamental order, subject to consistency with $tp_{\tTs}(a/Ae)$; this would
be even closer to the original proof.  

 \<{prop}\lbl{stableEI}  Let $T$ be a stable  theory admitting elimination of quantifiers, and let $\tTs$ be the theory
described above.     Assume $T$ eliminates imaginaries and,
and for any $A = \acl(A)$, $T_A$ eliminates
 finite generalized imaginaries.  Then  $\tTs$ admits elimination of imaginaries.  
\>{prop}
 
\<{remark}\lbl{stableEIr}  The stability condition can be weakened.  It suffices to assume $\Uu$ carries a notion of 
   independence with $2$-existence and uniqueness, and the following characterization of independence:
  (*) \ \    if $(a_i)$ is an indiscernible sequence over $A$, 
    $A_w=\acl(A\union \{a_i: i \in w \})$, and $A_w \meet A_{w'} = A$ for $w<w'$,  then the $(a_i)$
are independent over $A$.  
 \>{remark}

\proof    We may assume $T$ eliminates quantifiers.   By \propref{cond}, $T$ has 3-uniqueness.  
  It follows that $\tC=\{(A,\si): A \in \fC_T, \si \in Aut(A) \}$ has 3-existence (\propref{su}).
\def\tUu{\widetilde{\Uu}}
Let $\tUu$ be a qf-saturated model of $\tTs$.   

Part of the assumption is that finite sets are coded in $T$; hence also in $\tTs$.
Thus it suffices to prove elimination of imaginaries to the level of finite sets; in other words we have to show:   if $e$ is an imaginary element of $\tUu$, and $A$ is the set of real elements of $\acl_{\tTs} (e)$, then $e \in \dcl_{\tTs}(A)$.  

We have $e \in \dcl_{\tTs} (Aa)$ for some real tuple $a$; $e=a/E$ for some $\tUu$-definable equivalence relation $E$.  
  
Let $B$ be the set of real elements of $ \acl _{\tTs} (Aa)$.

If $b \in B \setminus A$, then $b \notin \acl(Ae)$, so $Aut(\Uu / Aeb)$ has infinite index in $Aut(\Uu/Ae)$.   It
follows that $\{g: g(b) =b\}$ is a subgroup of $Aut(\Uu/Ae)$ of infinite index. 
By  Neumann's Lemma
 \cite[Lemma 2.3]{neumann} , and compactness, it follows that there exists $g \in Aut(\Uu/Ae)$ such that $g(b) \neq b'$ 
for any $b,b' \in B \setminus A$. 

Let $a_1=a,a_2 = g(a)$, and define $a_n$ inductively so that 
$tp(a_n,a_{n+1}/A)=tp(a_1,a_2/A)$ and $a_{n+1}$  is independent 
from $a_1,\ldots, a_{n-1}$ over $A(a_n)$.  It is then easy to see that
$\acl(A,a_w) \meet \acl(A,a_{w'}) = A $ for any two   sets of indices $w,w'$ with $w<w'$.  
(Observe first that $\acl(A,a_w) \meet \acl(A,a_{w'}) \subset \acl(A,a_{m}) \meet \acl(A,a_{m+1})$, where
$m$ is the maximal element of $w$.)   However the $a_i$ need not be indiscernible; we convert them to an indiscernible sequence
in the paragraph below.  (Thanks to Christian d'Elb\'ee for pointing out that a previous treatment, simply citing Ramsey and compactness,  was inadequate.)


Let $u$ be a nonprincipal ultrafilter on $\Nn$,
and let $p$ be the $u$-average type of the sequence $a_n$.   Let $A'=\{a_1,a_2,\ldots\}$, and let $c_1 \models p| AA'$,
$c_2 \models p | AA'(c)$, etc.   Then the $(c_i)$ form an  indiscernible sequence, and $tp(c_i / A'(c_1,\ldots,c_{i-1}))$ is finitely satisfiable
in $A'$.  All the
 $c_n$ are $E$-equivalent since the $a_n$ are.       For any two   sets of indices $w,w'$ with $w<w'$, 
   $tp(c_{w'} / AA'(c_w)$ is finitely satisfiable in $A'$.   If  $d \in \acl(A,c_w) \meet \acl(AA') $ then $d \in \acl(A,a_v)$ for some finite $v \subset \Nn$;
   but then (internalizing $c_w$) we have  $d \in \acl(A,a_u)$ for some $u >v$; so $d \in A$.   
    If $d \in \acl(A,c_w) \meet \acl(A,c_{w'}) $ with $w < w'$ then  similarly we can see that $d \in \acl(A,a_v)$ for an appropriate $v$,
    and so $d \in  A$.   Thus (*) holds.  

By the assumption of  \remref{stableEIr}, 
$c_1, c_2$  are $A$- independent as tuples of $\Uu$, the $T$-restriction of $\tUu$.   
Hence by definition they are $A$-independent in $\tUu$.  
 We have found an $A$-independent pair of $E$-equivalent realizations of $tp(a/A)$.              
 
On the other hand, if $e \notin \dcl_{\tTs}(A)$, one easily obtains $E$-inequivalent 
independent elements realizing $tp(a/A)$.  
(E.g. $a$ and $g(a)$ where $g \in Aut_A(\tUu)$ and $g(e) \neq e$.)  Let $b \models tp(a/A)$ be such that $(c,c'), b$ are $A$-independent.
Then either $(c,b)$ or $(c',b) \notin E$.  

But a triangle with two equivalent and one inequivalent side cannot exist.  This contradicts 3-existence for $\tC$. 
\qed
%

\<{prop}\lbl{towards}   Let $T$ be a stable theory  admitting 
elimination of imaginaries.  Then $T$ has $4$-existence iff (with $A = \acl(\emptyset)$)
$T_A$   eliminates finite generalized imaginaries.  \>{prop}
   
\proof  
Since $T$ is stable, $T_A$ has $2$-uniqueness and hence $3$-existence over
any algebraically closed set $A'$.   By 
 \propref{cond}, $T$ has $3$-uniqueness; by \lemref{weakstepup},   it has $4$-existence.  
 
 Conversely, assume $T$ has
$4$-existence.  Let $G$ be a definable concrete groupoid with finite automorphism groups, defined in $T_A$.  
Fix a type $p$ of elements of $G$, and let $Tp$ be the set of types of
triples
$(a,b,d)$ with $a,b \models p$, $a,b$ independent over $A$,
and $c \in Mor_G(a,b)$.  Consider $(q_{12},q_{23},q_{13}) \in Tp$
such that there exist independent $a_1,a_2,a_3 \models p$
and $c_{ij}$ with $(a_i,a_j,c_{ij}) \models q_{ij}$ for $i<j$,
and such that $c_{12} = c_{23}^{-1}c_{13}$.  We can take $q_{23}=q_{13}$.
Such  triples  can be 4-amalgamated.  It follows easily that for
{\em any} independent $a_1,a_2,a_3$ and $c_{ij}$ with
$ (a_i,a_j,c_{ij}) \models q_{ij}$ for $i<j$, one has 
$c_{12} = c_{23}^{-1}c_{13}$.  (Otherwise, 3 triples with this
property and 1 triple without it could not be 4-amalgamated.)
Pick $(q_{12},q_{23},q_{13}) \in Tp$ with $q_{23}=q_{13}$.
It follows that for any independent $a_1,a_2 \models p$
there exists a {\em unique} $c \in Mor_G(a_1,a_2)$
with $q_{12}(a_1,a_2,c)$.  Moreover, we have a sub-groupoid
$G'$  of $G$ with the same objects and such that
$Mor_{G'}(a_1,a_2)$ is the unique realization of $q_{12}$.  
For any  functor $F$ on $G$ into definable sets, we now obtain
an equivalence relation on the disjoint union of the objects of $G$,
identifying $e \in F(a), e' \in F(a')$ if $h(e)=e'$ for the unique
$h \in \Mor_{G'}(a,a')$.  Using 
elimination of imaginaries, it is now easy to construct a finite group action equivalent to $G$.  \qed 

\<{cor} \lbl{fc} For stable $T$, the following are equivalent:  $4$-existence, $3$-uniqueness,
elimination of finite generalized imaginaries, $3$-existence for $\tC$.
  \>{cor}
\prf  By \propref{towards}, \propref{su}, \lemref{weakstepup}, \propref{cond}.  \eprf

\tsec{Discussion}  
 In many proofs regarding stable theories, there is no harm in passing to a theory $T'$ with more sorts, as long as $T$ remains stably embedded and with the same induced structure;
 especially if 
 $M' = \acl(M'|   {L}(T))$ for $M' \models T'$.   
  In this situation, by interpreting
 algebraic closure more widely in such extensions $T'$, the  $3$-uniqueness
or $4$-existence property for amalgamation holds.

A generalization of the above proof  for $n>3$ using an appropriate notion of higher groupoids, would be very  interesting.   \footnote{I have recently become aware of Jacob Lurie's work \cite{lure}, which may hold the key to this.
Note the apparent  resonance between Lurie's main theorem 6.1.0.6 there, and our \thmref{im-gr}.}

We do not  at present have a concrete description of the requisite sorts (analogous to equivalence relations or groupoids), but can at least prove their existence.

\<{prop}\lbl{higher}  Let $T$ be a theory with a canonical 2-amalgamation.  There exists an expansion $T^*$ of $T$ to a language with additional sorts, such that:

(1) $T$ is stably embedded in $T^*$, and the induced structure from $T^*$ on the $T$-sorts 
is the structure of $T$.  Each sort of $T^{*}$ admits
a 0-definable map to a sort of $T$, with finite fibers.

(2)  $T^*$ has existence and uniqueness for $n$-amalgamation. 
\>{prop}

We sketch the proof.

Condition (1) is equivalent to:  

(1')  If $N^* \models T^*$ and $N$ is the restriction to the sorts of $T$, then $Aut(N^*) \to Aut(N)$ is surjective, with profinite kernel.

For $T$ with canonical 2-amalgamation, and $p$  a type   of $T$ over $\acl(\emptyset)$, consider an expansion $T_p$ of $T$ as in \propref{cond}:
the points of a model $M_p$ of $T_p$ correspond to $\acl(a,M)$ where
$M$ is the restriction to $T$ of $M_p$, 
$a \models p$, and $a,M$ are embedded in some bigger model of $T$ via canonical 2-amalgamation.  

The proof of \propref{cond} shows that if each $T_p$ has unique
$(n-1,n+1)$-amalgamation then $T$ has unique $(n,n+1)$-amalgamation.

To prove the proposition, construct first an expansion $T^*$ of $T$ with property (1),
such that (U) any expansion of $T$ together with finitely many sorts of $T^*$ with property
(1) is equivalent to a sort of $T^*$.  

Note that $T^*_p$ enjoys the same property; since a relatively finite cover of $T^*_p$,
fibered over a sort $Y_a$, arises from a relatively finite cover of $T^*$, fibered over $Y$.

 Suppose   that unique $(n,n+1)$-amalgamation fails for a theory $T$ with the universal property
(U).  Take $n$ minimal.    Then   uniqueness at $(n-1,n+1)$ fails for $T_p$
for appropriate $p$.   
By \lemref{weakstepup}, $T'$ does not have $(n-1,n)$-uniqueness.    But $T_p$ also has (U).  
This contradicts the minimality of $n$.     

\<{problem}  For $T^*$, prove an analog of \secref{add-automorphism} for $n$ commuting automorphisms.  \>{problem}

 \>{section}

\<{section}{Linear   imaginaries} \lbl{linsec}
 
 We first discuss linear imaginaries in general; then restrict attention to the
triangular imaginaries that we will need.

\<{defn}  Let $\tres$ be a theory of fields (possibly with additional
structure.)  

 A  $\tres$-linear structure $\sA$ is a structure with a sort $\k$ for a model
 of $\tres$,  and additional sorts $V_i$ ($i \in I = I(\sA)$) 
 denoting finite-dimensional vector spaces.    Each $V_i$
has (at least) a $\k$-vector space structure, and $\dim_i V_i < \infty$.  
 
We assume:  \begin{enumerate}
  \item $\k$ is stably embedded, 

  \item the induced structure on $\k$ is precisely given by $\tres$
  \item  the $V_i$ are closed under tensor products and duals.

\end{enumerate}  \>{defn}

{\bf Explanation}  The language includes the language of $\tres$
(applying to $\k$), and for each $i$, a symbol for addition $+:V_i \times V_i \to V_i$, and scalar multplication $\cdot: k \times V_i \to V_i$.   

Given $i,j$, for some $\k$, the language includes  bilinear map
$b:V_i \times V_j \to V_k$, inducing an isomorphism $V_i \tensor V_j \to V_k$.

For each $i$, there is $j$ and a function symbol for
a pairing $V_i \times V_j \to\k$, inducing an isomorphism
$\dual{V_i} \to V_j$.

Additional
structure is permitted, subject to the embeddedness conditions (1,2).

Note that the tensor product $V \tensor W$ and dual $\dual{V}$ are at all events interpretable in $(k,V,W)$; so the conditions (3) can be
viewed as (partial) elimination of imaginaries conditions.

For any finite tuple $s=(s_1,\ldots,s_k)$ of indices, let 
$V_s = V_{s_1} \oplus \cdots \oplus V_{s_k}$, and let $P_s$
be the projectivization $V_s \- (0) / k^*$.  These can also clearly be viewed as imaginary sorts of $\sA$.

 \<{prop} \lbl{acflei}   Let $\k$ be an 
 algebraically closed field.  Then any $\k$-linear theory
 eliminates imaginaries to the level of the projective spaces
 $P_s$.    
\>{prop}

\proof This goes back to the 19th century (cf. references to Darboux in \cite{dar}) and occurs also in  
  \cite{acvf1}, Proposition 2.6.3, but not in easily quotable form.  If $W $ is a direct sum of some of the $V_i$, note first that the elements of the exterior powers $\Lambda^i W$ can be coded.  Indeed an element
of $\Lambda_i W$ can be viewed as a certain multilinear map on $\dual{W}$,
thus as an element of $W^{\tensor i}$.  This is again
a direct sum of some of the $V_i$.  

A $d$-dimensional subspace of such a $W$ corresponds to
a certain $1$-dimensional subspace of $\Lambda^k W$, and hence of
$W^{\tensor d}$.   Hence
it can always be coded as an element of the projectivization $P_s$.

Now any Zariski closed subset $Z$ of $W$ is determined, for some
$l$, by the space of polynomials of degree $\leq l$ vanishing on $Z$.
This is a subspace of $\bigoplus_{i \leq l} (\dual{W})^{\tensor i}$.  Hence
it is coded.  

It follows  by induction on dimension  that every definable subset of $W$ is coded (code the Zariski closure, and then the complement.)  \qed

A couple of remarks:

\<{lem} \lbl{elem} Let $A=(k,V_i)_{i \in I(A)}$ be a $\tres$ linear structure, and let
$k^*$ be an elementary extension of $\k$.  Let $V_i^* = k^* \tensor_k V_i$,
and let $A^* = (k^*,V_i^*)_{i \in I}$.  Then $A^*$ expands uniquely to an elementary extension of $A$.  (And every elementary extension of $A$ is obtained in this way.)  \>{lem}

\proof  Clearly, if $A \prec A^*$, then $\k \prec k^*$ and
(by the finite dimension)  $V_i^* = k^* \tensor_k V_i$.

Also, any $V_i$ has a basis $b_i$ in $A$.  There is a $b_i$-definable
bijection $f : V_i \to k^{n_i}$.  If  $R$ is a relation on $V_i$, or among several $V_i$, 
then $f  R$ is a relation $S$ on $\k$, and in any elementary extension
one must have:  $R = f^{-1}S$.   Thus uniqueness of the expansion is
clear, and it remains to show that this prescription always does yield
an elementary extension.   We may fix constants for each $b_i$.  But then
$V_i \subseteq dcl(k)$, and the assertion is immediate.  \qed

\<{lem} \lbl{constexp} Let $\tT$ be a theory, internal to a predicate $\k$, and with elimination of imaginaries. Let $\tT '$ 
be an expansion of $\tT$, such that every subset
of $k^m$, $\tT '$ -definable with parameters, is $\k$-definable with parameters.
Then $\tT '$ admits elimination of imaginaries.  \>{lem}

\proof    

\claim 
   If $X$
is a   definable subset of $\tT '$ (with parameters), then $X$ is 
also $\tT $ -definable  (with parameters).

\proof  By internality, there exists a $\tT$ definable
(with parameters) map $f$ on $\k^n$, whose image contains $X$.
$f^{-1}(X)$ is $\tT$-definable (with parameters) by the assumption
regarding new structure on $\k$.  Thus $X = ff^{-1}X$ is $\tT$-definable
(with parameters.)

 Hence any $\tT'$-definable set can be written $X=Y_c$ where $c$ is a canonical parameter
for $Y_c$ in $\tT$.  It follows that $Y_c$ is $\tT'$-definable and $c$ is
a canonical parameter
for $Y_c$ in $\tT '$. \qed

\tsec{Linear structures with flags}

We now consider {\em flagged spaces}.  For us this will mean:
a finite dimensional vector space $V$ together with a
 filtration $V_1 \subset \ldots \subset V_n =V $ by subspaces,
with $\dim V_i = i$. 

Given $V$, we form the dual $\dual{V}$ with the natural filtration
$\dual{V}_i$.  If $V$, $W$ are filtered spaces, take the tensor product with the filtration
$V_1 \tensor W_1 \subset V_1 \tensor W_2 \subset \cdots
 \subset V_1 \tensor W \subset  V_2 \tensor W_1 + V_1 \tensor W \subset \cdots \subset V_2 \tensor W \subset \cdots \subset V \tensor W$. 

Thus a family of flagged spaces can be closed under tensors and duals,
without losing the flag property.

\<{defn}  A  $\tres$-linear
structure $\sA$ {\em has flags} if:

(*)  For any $i$ with $\dim(V_i) >1$,  for some $j,k$ with $\dim(V_j)=\dim(V_i)-1$, $\dim V_k =1$, 
there exists a $0$-definable exact sequence
$0 \to V_k \to V_i \to V_j \to 0$.  \>{defn}

\<{lem} \lbl{projei}    (Elimination of projective imaginaries.) 
Let $\sA$ be a flagged $\tres$-linear structure.  Then 
elements of projectivizations of the vector spaces of $\tTt$
can be coded in $\tTt$.  In particular if $\k$ is an algebraically
closed field, $\tTt$   admits elimination of imaginaries.  \>{lem}

\proof  
 \lemref{acflei} applies here: 
 $(k,V_1 \subset \ldots \subset V_n)$
can be viewed as an expansion of $(k,V_1 \oplus \ldots \oplus V_n)$.
By \propref{acflei} and \lemref{constexp}, when $\k$ is algebraically closed,
all imaginaries of $\tTt$ are coded by elements of projective spaces.

So  it suffices to show in general how to code the projectivizations of the
vector spaces $W$ of $\tTt$.  Say 
$\dim(W)=d$.  $W$ 
 comes with   a   0-definable filtration, including 
$W_0 \subset W_1 \subset \cdots W_{\dim W} = W$ by
 subspaces $W_i$ with $\dim(W_i)=i$.   
Thus it suffices to code a 1-dimensional subspace $U$ of a filtered
vector space $W$.  Say $U \subset W_{k+1}, U \not \subseteq W_k$.  
Let $f: W_{k+1} \to Y = W_{k+1}/W_k  $ be the natural map.   Given $U$, one obtains $f|U$ and hence
$(f|U)^{-1}: Y \to W$.  But this is an {\em element} 
rather than a subspace of $Hom(Y,W) = \dual{Y} \tensor W$, and hence 
is coded.

\tsec{Linear structures with roots}

We say that a linear structure $A$ has roots if for any one-dimensional
$V=V_i$, and any $m \geq 2$, there exists $W = V_j$ and   $0$-definable
$\k$-linear embeddings
$f: W^{\tensor m} \to V_l$  and $g: V \to V_l$, with $g(V) \subset f(W)$.

A {\em good} linear structure is one with flags and roots.

\<{prop}  Let $\sA$ be a   linear structure with flags and roots, for an algebraically closed field $\k$
of characteristic $0$.
Then every finite internal cover of $\sA$ almost splits. \>{prop}

\proof    We may expand the theory by algebraic points; in particular we may assume
that the definable points of $\k$ form an algebraically closed field.
For $\k$ itself,
  the lemma follows from \remref{modelcov}.  For $\sA  $, we use \propref{finitecentral}.
  Let $D$ a $0$-definable set of $\sA$.
  $Aut(D/k)$ definably isomorphic to the $\k$-linear group $H$.  Since $\sA$ has
  flags, for each sort $W$ of $\sA$, $Aut(W/k) $ preserves a flag, so it is a solvable group, and hence so is $H$.
  Thus $H = UT$ where $U$ is the unipotent part of $H$.  Let $f: \tH \to H$ be a central
  extension of $H$ with prime cyclic kernel $Z= \Zz/l\Zz$.  We have to show that $\tH$ is
  represented as $Aut(\tD/k)$ for some $\tD$.
  
  We may take $H$ and $\tH$ connected.  Then $\tH$ is solvable.  So $\tH = \tilde{U} \tilde{T}$
  with $\tilde{U}$ unipotent.  $f | \tilde{U}$ has finite kernel, hence trivial kernel
  (  unipotent groups in characteristic zero have no nontrivial torsion elements.)    
 We obtain an induced map $\tilde{T} \to \tH / \tilde{U} \to H/U = T$,
 and see that $\tH \to H$ is induced from a central extension
 $0 \to Z \to \tilde{T} \to T \to 0$.  Now $T \iso (G_m)^n$,
 and every central extension of $T$ with kernel $Z$ is induced
 from the extension $G_m \to G_m$, $x \mapsto x^l$, via some
 rational character $\chi: T \to G_m$.   
 
 Let $D'$ be the set of bases (of the $\k$-space spanned by $D$)
 contained in $D$; then $Aut(D'/k)=H$ also; let $D'' = D' / Ker(\chi)$;
 by \lemref{acflei} and \lemref{projei}, $\sA$ admits elimination of 
 imaginaries, so $D''$ can be viewed as a 0-definable  set in $\sA$.
 Now $Aut(D''/k) = G_m$, and we have to represent the finite central 
 cover $x \mapsto x^l$.  
 
 $D''$ lies in some vector space $W$ in $\sA$.   Since $Aut(W/k) = G_m$,
 $W  = \oplus W_i$, and $G_m$ acts on $W_i$ via a character $\chi_i$.
 Each $1$-dimensional subspace $S$ of $W_i$ is $Aut(W/k)$-invariant.
  
 At the same time, $W$ contains a $0$-definable one-dimensional 
 subspace $W_1$; and $W/W_1$ is $0$-definable.  Either $S=W_1$,
 or $S$ embeds into $W/W_1$.  Continuing this way we find a $0$-definable
1-dimensional space $V$ such that $Aut(V/k)=Aut(S/k)$.  Doing
this for each $S$ in some some decomposition of $W$ into one-dimensional $Aut(W/k)$-invariant subspaces as above, we find $0$-definable
one-dimensional $V_1,\ldots,V_j$ with $Aut(W/k) = Aut(V_1,\ldots,V_j /k)$.

Now by taking roots of the $V_i$ we can find $\tV_i$ and
maps $(\tV_i)^{\tensor l} \to V_i$ such that
$Aut(\tV_i /k ) \to Aut(V_i /k)$ is isomorphic to $x \mapsto x^l,
G_m \to G_m$.  Pulling back to $Aut(\tV_1,\ldots,\tV_j /k)$ we succeed
in splitting the cover.   \qed
  
Note that a 1-dimensional vector-space $V$ amounts to the same thing
as a set $V^* = V \setminus (0)$ and a regular action of $k^*$
on $V^*$.  (Given such an action, recover the vector space structure on $V = V^* \union \{0\}$:
$0 \cdot v = 0$, $\alpha u + \beta u = (\alpha+\beta)u$.)  
So one has the usual $H^1$ formalism; the tensor product corresponds
to the sum in $H^1$.  

An $m$'th root $W$ of $V$, i.e. a one-dimensional vector space with
an isomorphism $f: W ^{\tensor m} \to V$, yields a  $(k^*)^m$
$T=T(W) \subset V^*$, i.e. a class of the equivalence relation of $(k^*)^m$-conjugacy.  
Namely,  let $t(w) = f(w \tensor \cdots \tensor w)$, and $T = t(W)$.  

Note that (because of $t$ and the $\k$- linear structure on $W$,
 $Aut(W/V,k) \leq \mu_m(k)$, the group of $m$'th roots of unity in $\k$.  Thus $W$ can be regarded as 
 a finite internal cover of $(V,k)$, and in particular as a generalized
 imaginary sort.  
 

\<{remark} \lbl{leirem} Assume $\tres$ admits linear elimination of imaginaries for structures with flags and roots,
and let $\sA$ be a $\tres$-linear structure with flags.  Assume that
for any 1-dimensional $V$ of $\sA$, and any $m$, $V$ has
a distinguished $(k^*)^m$-class (over $\emptyset$.)  Then $\sA$
admits elimination of imaginaries. \>{remark}

\proof  We can arrive at a structure with EI by successively adding roots
to one -dimensional vector spaces $V$ of $\sA$.  This involves
adding an $m$'th root $W$, all tensor powers $W^{\tensor n}$
for $n \in \Zz$, and all tensor products $W^{\tensor n} \tensor U$
for $U$ a vector space of $\sA$.   If 
$\mu_m(k)=1$, then $W$ can be identified with a $(k^*)^m$-class
in $V$, and $W^{\tensor n} \tensor U$ embeds into $V^{\tensor n} \tensor U$ via this class.  Otherwise, $\mu_m(k)$ acts on the new structure,
fixing $\k$ and $\sA$, 
and one sees that if an element $w \tensor u \in W^{\tensor n} \tensor U$ is fixed by $\mu_m(k)$, then $n=0$, or $w=0$.  Thus any imaginary
of $\sA$ coded in the new structure already lies in $\sA$.   Applying
this iteratively, we see that $\sA$ has EI.   \qed
 
\
\ssec{Pseudo-finite fields}  \lbl{pseudofin}

 A {\em pseudo-finite} field is a perfect PAC field (i.e. every irreducible variety over
 $F$ has an $F$-point), with  Galois group $\hat{\Zz}$.  Ax showed
 that $F$ is pseudo-finite iff every sentence true in $F$ is true in infinitely many finite fields.
 (\cite{ax}, \cite{FJ}, \cite{pac}.)
 {\em We take $F$ to come together with an isomorphism $\hat{\Zz} \to Gal(F)$.}  
 In terms of language, this means that we fix, for each $d \in \Nn$, an imaginary element coding
 a generator of $Gal(F_d / F)$ (where $F_d$ is the extension of $F$ of order $d$.)
 (This  is a little more than the pure field structure; but  all finite fields 
 do have a canonical generator of Galois, so perhaps  it's only fair that the pseudo-finite ones should.)
 When $char(k)=0$, this
   is equivalent to fixing a surjective group homomorphism $ k[\mu_d]^* \to \mu_d$;
 cf. \cite{CH}.   As noted there, in this language, $F$ admits elimination of
 imaginaries.   \footnote{It is also asserted there that $F$ admits elimination
 of imaginaries in the field language, up to sorts coding elements of the Galois
 group; but  this slighlty stronger statement is incorrect.}

\<{thm} \lbl{pfeli} Let  $F$ be a pseudo-finite field of characteristic $0$ (with fixed generator of Galois).    Then a good linear structure over $F$
admits   elimination of imaginaries.  \>{thm}

\proof   Let $\tres = Th(F)$, and let $\sA$ be a good $\tres$-linear
structure.   If $\sA ' $ is a reduct of $\sA$, with the same sorts, with the full structure on
the field sort, and remembering the $F$-linear structure of each vector space in $\sA$
and all $0$-definable linear maps among them, then $\sA'$ is also good.  
Moreover $\sA$ admits EI if $\sA'$ does  (\lemref{constexp}.)  Thus we may assume
the structure on $\sA$ consists just of the structure on $\k$, the $\k$-linear structures
and the $0$-definable $\k$-linear maps among them.  
 
Let $\sA ^a$ be the linear structure over the algebraic closure 
$F^a$, and with vector spaces $V^a = F^a \tensor_F V$ for each vector space
$V$ of $\sA$.  Any $0$-definable linear map in $\sA$ extends uniquely to a
$0$-definable linear map in $\sA ^a$, and we take this to define a structure on $\sA ^a$.  
Since $\sA$ is good, so is $\sA^a$.    Let $T_1 = Th(F^a, V^a)_{V \in \sA}$.

 \begin{enumerate}
 
  \item   $T_1$ admits elimination of imaginaries 
  (\lemref{acflei}) and of generalized finite imaginaries,
  and is stable (every sort has finite Morley rank.)

  \item  Let $T_2$ be the model companion of the  theory of pairs $(M,\si)$ where $M \models T_1$ and $\si \in Aut(M)$.  Then $T_2$ admits elimination of imaginaries.   (\lemref{stableEI}.)

\item  Let $(K,V,\si) \models T_2$.  Let $(F,V)$ be the fixed field and fixed vector spaces
of $\si$.  Then $\dim_F V = d$.   $F$ is pseudo-finite, and $K$ 
 can be chosen so that the fixed field will have the same theory  as  the original $F$.  
 (cf.  \cite{CH}.)
 
 \item  Any imaginary of $(F,V)$ is in particular an imaginary of $(K,V,\si)$, and thus
 can be coded by a tuple of elements of 
 $(K,V)$.  Each of these elements must
 be fixed by $\si$.
 
 \item  $(F,V_F)$ coincides with the $\si$-fixed part of $(K,V)$.   This follows
 from the fact that $V$ has a $\si$-fixed basis.
 
 \end{enumerate}
  \qed.
          \>{section}


\begin{thebibliography}{99}
 
 \bibitem{ah-z} Ahlbrandt, Gisela; Ziegler, Martin Quasi-finitely axiomatizable totally categorical theories. Stability in model theory (Trento, 1984). Ann. Pure Appl. Logic 30 (1986), no. 1, 63--82. 
 
\bibitem{ah-z1}  Ahlbrandt, Gisela; Ziegler, Martin What's so special about $({\bf Z}/4{\bf Z})^\omega$(Z/4Z)?? Arch. Math. Logic 31 (1991), no. 2, 115Ð132. 
 
  \bibitem{ax}  
 Ax, James The elementary theory of finite fields.  Ann. of Math. (2) 88 1968 239--271.
 

 
\bibitem{baldwin}  Baldwin, John T.  \textsl{Fundamentals of stability theory.}  Perspectives in Mathematical Logic. Springer-Verlag, Berlin, 1988. xiv+447 pp. ISBN: 3-540-15298-9 
  
\bibitem{bouh} ÊBouscaren, E.; Hrushovski, E.,  On one-based theories.  J. Symbolic Logic 59 (1994), no. 2, 579--595.
  
 \bibitem{CH}  
Chatzidakis, Z.,
Hrushovski, E., Model Theory of difference fields, AMS Transactions v.
   351, No. 8, pp. 2997-3071 

 \bibitem{CP}   Chatzidakis, Z., Pillay, A.,
Generic structures and simple theories.   
Ann. Pure Appl. Logic 95 (1998), no. 1-3, 71--92

 \bibitem{ch-h} Cherlin, Gregory; Hrushovski, Ehud,Finite structures with few types. 
Annals of Mathematics Studies, 152. 
Princeton University Press, Princeton, NJ, 2003. vi+193 pp. ISBN 0-691-11331-9
 

\bibitem{evh}  Evans, David M.; Hrushovski, Ehud On the automorphism groups of finite covers. Stability in model theory, III (Trento, 1991).  Ann. Pure Appl. Logic  62  (1993),  no. 2, 83--112.


\bibitem{ev}  Evans, David M. Finite covers with finite kernels. Joint AILA-KGS Model Theory Meeting (Florence, 1995).  Ann. Pure Appl. Logic  88  (1997),  no. 2-3, 109--147. 


\bibitem{emi}   Evans, David M.; Macpherson, Dugald; Ivanov, Alexandre A. Finite covers.  Model theory of groups and automorphism groups (Blaubeuren, 1995),  1--72, London Math. Soc. Lecture Note Ser., 244, Cambridge Univ. Press, Cambridge, 1997. 
 

\bibitem{FJ}
\textsc{M. Fried and M. Jarden:} \textsl{Field
Arithmetic}, Springer Verlag, Berlin
1986.


 \bibitem{acvf1} D. Haskell, E. Hrushovski, H.D. Macpherson, `Definable sets in 
algebraically closed valued fields. Part I: elimination of imaginaries', submitted.

 \bibitem{acvf2} D. Haskell, E. Hrushovski, H.D. Macpherson, `Definable sets in 
algebraically closed valued fields,
 Part II: stable domination and independence.  Submitted.  
 
 
 

 
 \bibitem{pac}    Hrushovski, Ehud Pseudo-finite fields and related structures.  Model theory and applications,  151--212, Quad. Mat., 11, Aracne, Rome, 2002
 
 \bibitem{bedlewo}  Ehud Hrushovski, Computing the Galois group of a linear differential equation, 
in {\em Differential Galois Theory}, 
Banach Center Publications 58, Institute of Mathematics, 
Polish Academy of Sciences, Warszawa 2002

 \bibitem{dar} 
E. Vessiot, M\'ethodes d'integration \'elementaires, {\it in} {\bf Encyclop\'edie
des Sciences Math\'ematiques Pures et Appliqu\'ees}, ed. Jules Molk, Tome II (3\'eme vol), 
Gauthier-Villars Paris / B.G. Teubner , Leipzig, 1910; 
reprinted by \'Editions Jacques Gabay, 1992


 \bibitem{HP} Hrushovski, E. \& Pillay, A., "Groups Definable in Local
Fields and Pseudo-Finite Fields", {\em Israel Journal of Mathematics}
{\bf85} (1994) pp. 203-262

  
 
 
 \bibitem{meta} Hrushovski, E., valued fields, metastable groups (preprint)
 
  \bibitem{kpy}   Byunghan Kim, T. de Piro and J. Young, {\em Constructing the hyperdefinable group from the group configuration}, http://www-math.mit.edu/~bkim/
 
\bibitem{lurie} Lurie, Jacob Higher topos theory. Annals of Mathematics Studies, 170. Princeton University Press, Princeton, NJ, 2009. xviii+925 pp.  See also {\text http://www.math.harvard.edu/~lurie/}.


\bibitem{neumann} P.M. Neumann, `The structure of finitary permutation groups',
 Archiv der Mathematik 27 (1976), 3--17.

 
\bibitem{poizat} \textsc{B. Poizat:} Une th\'eorie de Galois
imaginaire, \textsl{J. Symbolic Logic} \textbf{48}
(1983) 1151-1170.

\bibitem{lazy} Shelah, Saharon
The lazy model-theoretician's guide to stability.
Comptes Rendus de la Semaine d'\'Etude en Th\'eorie des Mod\`eles (Inst. Math., Univ. Catholique Louvain, Louvain-la-Neuve, 1975).
Logique et Analyse (N.S.) 18 (1975), no. 71-72, 241--308.

\bibitem{shelah} \textsc{S. Shelah:}  \textsl{Classification Theory and the number of non-isomorphic models}, revised edition, North-Holland Amsterdam-Tokyo 1990.

\bibitem{sh87b}  \textsc{S. Shelah:}Classification theory for nonelementary classes, I. The number of uncountable models of $\psi \in L_{\omega _{1},\omega }$, Part B. 
 \textsl{Israel J Math} {\bf 46} (1983) 241-273
 

 
 \bibitem{wagner} Wagner, Frank O. Simple theories. Mathematics and its Applications, 503. Kluwer Academic Publishers, Dordrecht, 2000.
 
\end{thebibliography}
\end{document}